\documentclass[12pt,a4paper]{amsart}
\usepackage{amsmath}
\usepackage{mathabx}
\usepackage{amstext}
\usepackage{amssymb}
\usepackage{array}
\usepackage{amsthm}
\usepackage{dsfont}
\usepackage[latin1]{inputenc}
\usepackage[T1]{fontenc}
\usepackage{mathtools}
\usepackage{appendix}
\usepackage{xcolor}
\usepackage{ulem}
\usepackage[arrow, matrix, curve]{xy}
\usepackage{verbatim}
\usepackage{graphicx}
\usepackage{subfigure}
\usepackage{hyperref}
\usepackage{cases}
\usepackage[update,prepend]{epstopdf}
\usepackage{lmodern}

\newtheorem{theorem}{Theorem}
\theoremstyle{plain}
\newtheorem{definition}[theorem]{Definition}
\newtheorem{assumption}[theorem]{Annahme}
\newtheorem{lemma}[theorem]{Lemma}

\newtheorem{proposition}[theorem]{Proposition}
\theoremstyle{definition}

\numberwithin{equation}{section}

\def\eps{\varepsilon}

\def\ri{{\rm i}}

\def\uhat{\widehat{u}}
\def\util{\tilde{u}}
\def\utilext{\tilde{u}_\text{ext}}
\def\ubartilext{\tilde{\widebar{u}}_\text{ext}}
\def\vtil{\tilde{v}}

\def\Ahat{\widehat{A}}
\def\Ahatbar{\widehat{\widebar{A}}}
\def\Uvec{\vec{U}}
\def\Uvecext{\vec{U}^\text{ext}}
\def\Uvecapp{\vec{U}^\text{app}}
\def\Uapp{U^\text{app}}
\def\uapp{U_\text{app}}
\def\Uext{U^\text{ext}}
\def\Hsper{H^s_\text{per}}

\def\bi{\begin{itemize}}
\def\ei{\end{itemize}}

\newcommand{\R}{\mathbb{R}}
\newcommand{\C}{\mathbb{C}}
\newcommand{\B}{\mathbb{B}}
\renewcommand{\P}{{\mathbb P}}

\newcommand{\N}{\mathbb{N}}

\newcommand{\Z}{\mathbb{Z}}

\newcommand{\cD}{{\mathcal D}}
\newcommand{\cE}{{\mathcal E}}

\newcommand{\cL}{{\mathcal L}}

\newcommand{\cT}{{\mathcal T}}

\newcommand{\cX}{{\mathcal X}}

\renewcommand{\phi}{\varphi}

\newcommand{\pa}{\partial}

 \def\dd{\, {\rm d}}

\DeclareMathOperator{\sign}{sign}

\DeclareOldFontCommand{\it}{\normalfont\itshape}{\mathit}

\newcommand{\bspm}{\left(\begin{smallmatrix}}\newcommand{\espm}{\end{smallmatrix}\right)}
\newcommand{\bpm}{\begin{pmatrix}}\newcommand{\epm}{\end{pmatrix}}

\def\bs{\begin{satz}}\def\es{\end{satz}}
\def\blem{\begin{lemma}}\def\elem{\end{lemma}}
\def\bthm{\begin{theorem}}\def\ethm{\end{theorem}}
\def\bcor{\begin{corollary}}\def\ecor{\end{corollary}}
\def\beq{\begin{equation}}\def\eeq{\end{equation}}
\def\beqq{\begin{equation*}}\def\eeqq{\end{equation*}}
\def\bal{\begin{align}}\def\eal{\end{align}}
\def\bpf{\begin{proof}}\def\epf{\end{proof}}
\def\bex{\begin{example}}\def\eex{\end{example}}
\def\brem{\begin{remark}}\def\erem{\end{remark}}
\def\bass{\begin{assumption}}\def\eass{\end{assumption}}
\def\bprop{\begin{proposition}}\def\eprop{\end{proposition}}
\def\bdefi{\begin{definition}}\def\edefi{\end{definition}}

\DeclareFontEncoding{FMS}{}{}
\DeclareFontSubstitution{FMS}{futm}{m}{n}
\DeclareFontEncoding{FMX}{}{}
\DeclareFontSubstitution{FMX}{futm}{m}{n}
\DeclareSymbolFont{fouriersymbols}{FMS}{futm}{m}{n}
\DeclareSymbolFont{fourierlargesymbols}{FMX}{futm}{m}{n}
\DeclareMathDelimiter{\VERT}{\mathord}{fouriersymbols}{152}{fourierlargesymbols}{147}

\def\bi{\begin{itemize}}
\def\ei{\end{itemize}}
\def\ben{\begin{enumerate}}
\def\een{\end{enumerate}}
\usepackage{enumitem}
\setitemize{leftmargin=*}

\begin{document}

\title[NLS Approximation for Wavepackets in $\R^d$]{NLS approximation for wavepackets in periodic cubically nonlinear wave problems in $\R^d$}
\author{Tom\'a\v{s} Dohnal$^{1}$ and Daniel Rudolf$^{2}$}
\address{$^{1}$ Technische Universit\"at Dortmund, Fakult\"at
      f\"ur Mathematik, Vogelpothsweg 87, D-44227 Dortmund, Germany; $^{2}$ Ruhr-Universit\"at Bochum, Fakult\"at
      f\"ur Mathematik, Universit\"atsstrasse 150, D-44801 Bochum, Germany  
      Germany.}
\email{tomas.dohnal@math.tu-dortmund.de, daniel.rudolf@ruhr-uni-bochum.de}
\date{\today}


\begin{abstract}
The dynamics of single carrier wavepackets in nonlinear wave problems over periodic structures can be often formally approximated by the constant coefficient nonlinear Schr\"odinger equation (NLS) as an effective model for the wavepacket envelope. We provide a detailed proof of this approximation result for the Gross-Pitaevskii equation (GP) and a semilinear wave equation, both with periodic coefficients in $\N\ni d$ spatial dimensions and with cubic nonlinearities. The proof is carried out in Bloch expansion variables with estimates in an $L^1$-type norm, which translates to an estimate of the supremum norm of the error. The regularity required from the periodic coefficients in order to ensure a small residual and a small error is discussed. We also present a numerical example in two spatial dimensions confirming the approximation result and presenting an approximate traveling solitary wave in the GP with periodic coefficients.
\end{abstract}

\maketitle

   \vspace*{2mm} {\bf Key-words:} wavepackets, NLS approximation, Bloch waves, periodic media, solitary wave, approximation error, Gross-Pitaevskii, nonlinear wave equation

     \vspace*{2mm}
		
\setcounter{page}{1}

\section{Introduction}
We study the asymptotics of wavepackets in $\N\ni d$-dimensional cubically nonlinear wave problems with periodic coefficients. We consider wavepackets given by a single carrier wave modulated by a slowly varying envelope. In periodic media the carrier wave is a Bloch wave of the corresponding linear problem. For envelopes which are also appropriately small in amplitude an effective constant coefficient equation of the nonlinear Schr\"odinger (NLS) type can be easily derived. This equation depends only on the slow variables of the envelope and not on the fine oscillations in the carrier wave. In order to guarantee that solutions of the effective equation produce an approximation of solutions of the original problem, the approximation error needs to be estimated.

The motivation for this analysis is twofold. Firstly, it provides a rigorous justification of the approximation of wavepackets in nonlinear periodic media in arbitrary dimension via a simple constant coefficient equation in the slow variables. The fact that only slow variables appear means also that much coarser discretizations can be used when solving the problem numerically. Secondly, because the effective NLS-type equation typically supports solitary wave solutions, the results produce approximate moving solitary waves in periodic media. The velocity of the solution is asymptotically close to the group velocity of the carrier Bloch wave. Moreover, certain (e.g. radially symmetric) solitary waves of the effective NLS equation can be easily found numerically. By choosing carrier waves with different group velocities, we can tune the velocity of propagation of the pulse. Hence, in principle, close to solitary pulses are produced in $d$-dimensional periodic media for any propagation direction. Stability of such solutions is, of course, of importance but outside the scope of this paper.

We consider two standard models: the Gross-Pitaevskii equation (GP) 
$$\ri \pa_t u +\Delta u - V(x) u -\sigma(x) |u|^2 u =0, \quad (x,t)\in \R^d\times \R$$
and the  semilinear wave equation (NLW) 
$$
\pa_t^2 u=\chi_1(x)\Delta u-\chi_2(x)u-\chi_3(x)u^3,  \quad (x,t)\in \R^d\times \R.
$$
The coefficients $V,\sigma,\chi_j$ for $j=1,2,3$ are chosen $2\pi-$periodic in each coordinate for simplicity. 

The Gross-Pitaevskii equation (sometimes called the periodic nonlinear Schr\"odinger equation) describes Bose-Einstein condensates superimposed onto an optical lattice
\cite{Gross61,Pitaev61,OK03} but it is also an approximate model for light propagating in Kerr-nonlinear photonic crystals \cite{MK01,SK02,Efre03,KEVT06}. 

The semilinear wave equation is a simplified model corresponding to a reduction of the nonlinear Maxwell equations for Kerr-nonlinear photonic crystals. The main difference is that the Maxwell equations are quasilinear since they include the term $\pa_t^2(u^3)$. The above NLW has been considered in the same context, for instance, in \cite{BSTU06}. A proof of the approximation result for a quasilinear equation with $\pa_t^2(u^3)$ has been presented in \cite{LS12}.

Solitary waves, close-to-solitary waves or, more generally, localized coherent waves are of phenomenological as well as applied interest. Mainly in telecommunication and optical computing such pulses are promising as possible bit-carriers. Periodic media in optics, so called photonic crystals, are being considered as components in future optical computing, switching and logic.  
Nonlinear pulses have been observed in photonic crystals \cite{SE03} and have been proposed as bit carriers in the above applications \cite{Brod_98,Pelusi_etal_08}. Also in Bose-Einstein condensates solitary waves have been observed \cite{Eiermann_etal_04}.

The mathematical question of whether an effective equation for the envelope of a wavepacket truly approximates the dynamics is not trivial. A formally derived effective model can indeed fail to approximate a solution of the original problem \cite{Schn1995,GS01,SSZ2015}. We show here that for GP and NLW an approximation holds under some regularity assumptions, a simpleness assumption on the eigenvalue of the corresponding Bloch eigenvalue corresponding to the carrier and in the NLW case under a spectral non-resonance condition. We approximate wavepackets of GP via the ansatz
\beq\label{E:ans}
u(x,t)\approx \eps A(\eps(x-v_g t),\eps^2t)p_{n_0}(x,k_0)e^{\ri(k_0\cdot x-\omega_0 t)}
\eeq
and via two times the real part of the above for the NLW. Here $\eps>0$ is a small asymptotic parameter, $A$ is the (unknown) envelope, $p_{n_0}(x,k_0)e^{\ri(k_0\cdot x-\omega_0 t)}$ is a carrier Bloch-wave, and $v_g$ is its group velocity. Hereby $(\omega_0, p_{n_0})$ is an eigenpair of the Bloch eigenvalue problem and $k_0$ is a wavevector in the Brillouin zone; for details see below. The dependence of $A$ on $\eps^2 t$ models a slow temporal modulation of the wavepacket resulting from a periodicity induced dispersion and from the nonlinearity. The result says, roughly speaking, that if the envelope $A$ satisfies the effective NLS equation
$$\ri \pa_TA(X,T) +\frac{1}{2}\nabla \cdot (D^2\omega_{n_0}(k_0)\nabla A(X,T)) +\nu (|A|^2A)(X,T)=0, \quad (X,T) \in \R^d\times \R,$$
then for all $\eps$ small enough the above ansatz approximates a true solution $u$ on time intervals $[0, c\eps^{-2}]$. Here $X=\eps(x-v_gt), T=\eps^2t$, $k\mapsto (\omega_{n}(k))_{n\in \N}$ is the band structure, and $\nu\in \R$ is an integral of the Bloch eigenfunction $p_{n_0}(\cdot,k_0)$. The precise statements of our results are Theorems \ref{S:GP} and \ref{S:NLW}. In one spatial dimension the approximation has been proved for NLW in \cite{BSTU06}. In \cite{Pelinov_2011} the (technically simpler) GP case was considered - but only for $v_g=0$.

Other wavepackets than those given by ansatz \eqref{E:ans} are possible. In particular, one can use an ansatz with several carrier Bloch waves each modulated by an envelope. Because the carrier waves typically have different group velocities, the envelopes do not depend on one moving frame variable, but rather simply on $\eps x$ and $\eps t$. The resulting effective equations are of first order, so called coupled mode equations (CME). In the GP setting the most general ansatz was considered in \cite{GMS08}, where, however, the question of solitary waves was not considered. In one dimension CME possess a spectral gap and families of solitary waves, so called gap solitons, exist \cite{AW89}. Earlier, CME have been rigorously justified in one dimensional periodic structures with infinitesimal contrast in \cite{GWH01,SU01,Pelinov_2011}. In \cite{DH17} the proof for the case of arbitrary contrast in one dimension was provided in a setting leading to gap solitons.

Our approach to justifying the NLS asymptotics for the ansatz \eqref{E:ans} is similar to that in \cite{BSTU06} for the one dimensional case but we provide a more detailed analysis of the residual. We also discuss the required regularity of the periodic coefficients $V,\sigma$, resp. $\chi_j,j=1,2,3$. The required regularity grows with the dimension $d$.
The proof uses first the Bloch transformation (a generalization of the Fourier-transform), which transforms the problem from the domain $\R^d$ to a $d-$dimensional torus with a wavenumber parameter $k$ in the Brillouin zone $\B$. Next, an expansion of the solution in the eigenfunctions of the Bloch eigenvalue problem is applied, which transforms the problem to a system of ordinary differential equations in time parametrized by $k\in\B$. As the Bloch transform in space is an isomorphism between $H^s(\R^d)$ and $L^2(\B,H_\text{per}^s(\P))$ and the expansion in Bloch eigenfunctions is an isomorphism between $L^2(\B,H_\text{per}^s(\P))$ and $L^2(\B,l^2_{s/d})$, the authors of \cite{BSTU06} work in $L^2(\B,l^2_{s/d})$ in the Bloch-expansion variables and produce thus an error estimate in $H^s(\R^d)$. In higher dimensions, however, estimates in $L^2$ (and hence also in $H^s$) are more difficult due to the loss of $\eps$ powers when evaluating the norm of terms of the form $f(\eps x,t)$. It is namely $\|f(\eps \cdot,t)\|_{L^2(\R^d)}=\eps^{-d/2}\|f(\cdot,t)\|_{L^2(\R^d)}$. This has the effect that for satisfactory $H^s$-estimates of the error higher order correction terms have to be included in the ansatz in order to produce a residual that is small enough as $\eps\to 0$. To avoid this, we work in $L^1(\B,l^2_{s/d})$ instead of $L^2(\B,l^2_{s/d})$. Although no obvious isomorphism holds between this $L^1$-space and a space in the physical variables, the supremum $x-$norm is controlled by the $L^1(\B,l^2_{s/d})$-norm. The use of the $L^1$ space for the justification of amplitude equations appears in \cite{SU01,PS07} as well as, for example, in Chapter 2 of \cite{Pelinov_2011}.

The approaches for the equations GP and NLW are analogous but while the nonlinearity $|u|^2u$ in the GP case is gauge invariant, in the NLW case the nonlinearity $u^3$ applied to the ansatz \eqref{E:ans} generates higher harmonics, which have to be accounted for in a refined ansatz. Eliminating the leading order part of the residual at these higher frequencies also leads to a non-resonance condition for the linear spectral problem. This is an implicit condition on the linear coefficients $\chi_1$ and $\chi_2$.

The rest of the paper is structured as follows.  In Sec. \ref{S:FA} we collect the functional analytic results needed in the proof, in particular we review the Bloch transformation, the Bloch eigenvalue problem and the expansion in Bloch eigenfunctions and discuss the regularity of Bloch eigenfunctions. In Sec. \ref{S:GP} we prove the approximation result for the Gross-Pitaevskii (GP) equation and in Sec. \ref{S:NLW} for the nonlinear wave equation (NLW). In each case we define a modified (extended) approximation ansatz, for which we calculate the residual and derive the effective NLS equation. Next, we estimate the residual and the approximation error. Finally, Sec. \ref{S:num} presents a numerical example for the GP in two dimensions, confirming the $\eps$-convergence of the error using a solitary wave solution as an example. 

\section{Functional analytic tools and lemmas}\label{S:FA}

\subsection{Bloch transformation}
We recall first the definition and some basic properties of the Bloch transformation. For further information see e.g. \cite{PBL-1978}.

Using the Fourier transformation
$$f\mapsto \hat{f}, \ \hat{f}(k)=\frac{1}{(2\pi)^d}\int_{\R^d}f(x)e^{-\ri k\cdot x}\dd x$$
with the inverse $f(x)=\int_{\R^d}\hat{f}(k)e^{\ri k\cdot x}\dd k$, one can formally write
\beqq f(x)=\int_{\R^d} \widehat{f}(k)e^{\ri k\cdot x}\mathrm{d}k=\sum_{m\in\Z^d}\int_\B \widehat{f}\left(k+m\right)e^{\ri\left(k+m\right)\cdot x}\mathrm{d}k=\int_\B \widetilde{f}(x,k)e^{\ri k\cdot x}\dd k\eeqq
with $\B:=\left(-\tfrac{1}{2},\tfrac{1}{2}\right]^d,$ and
\beqq \widetilde{f}(x,k):=\sum_{m\in\Z^d}\widehat{f}\left(k+m\right)e^{\ri m\cdot x}.\eeqq
The Bloch transformation defined for $s\geq 0$ by                   
\beqq \cT:H^s(\R^d)\rightarrow L^2(\B,\Hsper(\P)), \quad f\mapsto \cT (f):=\widetilde{f}\eeqq
is an isomorphism \cite{RS4} with the inverse $f(x)=\int_\B \widetilde{f}(x,k)e^{\ri k\cdot x}\dd k$. Here 
$\P:=(0,2\pi]^d$ and
$\Hsper(\P)$ is the closure of $\P$-periodic $C^\infty(\R^d)$ functions in the $H^s(\P)$-norm.
A direct calculation shows that the following properties hold for all $f,g\in H^s(\R^d)$ and $1\leq j \leq d$, $\N \ni p \leq s$:
\beq\label{E:perinxandk}
\cT (f)(x+2\pi e_j,k)=\cT (u)(x,k),\quad \cT (f)\left(x,k+e_j\right)=e^{-\ri x_j}\cT(f)(x,k),
\eeq
\beq\label{E:xandtder}
\cT (\partial_{x_j}^pf)(x,k)=(\partial_{x_j}+\ri k_j)^p(\cT f)(x,k),
\eeq
\beq\label{E:commper}
\cT(Vf)(x,k)=V(x)\cT(f)(x,k), \text{ if } V(x)=V(x+2\pi e_j)\text{ for }1\leq j\leq d,
\eeq
\beq\label{E:conv}
\cT(fg)(x,k)=((\cT f)*_{\B} (\cT g))(x,k):=\int_{\B}(\cT f)(x,k-l)(\cT g)(x,l)\mathrm{d}l.
\eeq

In our analysis, however, we do not use $L^2$-estimates of the residual and the asymptotic error since too many $\eps$-powers are lost in $L^2$. This can be seen from the fact that the residual consists of functions of the form $f(\eps x)g(x)$ and $\|f(\eps \cdot)g(\cdot)\|_{L^2(\R^d)}\leq \eps^{-d/2}\|g\|_{L^\infty(\R^d)}\|f\|_{L^2(\R^d)}$. This loss of $\eps$-powers means that without including higher order correction terms in the ansatz the resulting asymptotic error is not $\mathcal{O}(\eps^2)$ in $\|\cdot\|_{H^s}$ on the desired time interval $[0,c\eps^{-2}]$.

Instead, we work in $L^1$ in the Bloch variable $k$, which results in supremum-norm estimates in the physical $x$-variables. In detail we use the norm
\beq \|\widetilde{u}\|_{L^1(\B,\Hsper(\P))}:=\int_{\B}\|\widetilde{u}(\cdot,k)\|_{H^s(\P)}\mathrm{d}k.\label{DefL1}\eeq
Although we lose the above isomorphism-property, Lemma \ref{L:sup_control}
guarantees the control of the supremum norm of $u$ via $\|\widetilde{u}\|_{L^1(\B,\Hsper(\P))}$.
\blem\label{L:sup_control}
Let $s>d/2$.  There is $c>0$ such that for all $\util \in L^1(\B,\Hsper(\P))$ which satisfy \eqref{E:perinxandk}, we have for the function $u(x):=\int_\B \util(x,k)e^{\ri k\cdot x}\dd k$
$$\sup_{x\in \R^d}|u(x)|\leq c\|\util\|_{L^1(\B,\Hsper(\P))} \quad \text{and} \ u(x)\to 0 \text{ as } |x| \to \infty.$$
\elem
The proof is completely analogous to that of Lemma 2 in \cite{DH17}.

As the next lemma shows, the algebra property of $H^s$ yields the algebra property also for the space $L^1(\B,\Hsper(\P))$.
\blem\label{L:algeb_L1Hs}
Let $\util,\vtil\in L^1(\B,\Hsper(\P))$ with $s>d/2$. Then
$$\|\util \ast_{\B}\vtil\|_{L^1(\B,\Hsper(\P))} \leq c \|\util\|_{L^1(\B,\Hsper(\P))}\|\vtil\|_{L^1(\B,\Hsper(\P))}.$$
\elem
\bpf
$$
\begin{aligned}
\|\util \ast_{\B}\vtil\|_{L^1(\B,\Hsper(\P))} & \leq c\int_{\B}\int_{\B}\|\util(\cdot,k-l)\|_{H^s(\P)}\|\vtil(\cdot,l)\|_{H^s(\P)} dldk \\
&\leq c \|\util\|_{L^1(\B,\Hsper(\P))}\|\vtil\|_{L^1(\B,\Hsper(\P))},
\end{aligned}
$$
where the first inequality follows by the algebra property of $H^s$ in $d$ dimensions, i.e. $\|fg\|_{H^s(\Omega)}\leq c\|f\|_{H^s(\Omega)}\|g\|_{H^s(\Omega)}$ for $\Omega\subset \R^d$ and $s>d/2$, see Theorem 5.23 in \cite{Adams}. The second step follows by Young's inequality for convolutions.
\epf

\subsection{Bloch eigenvalue problem}

Our two examples (GP and NLW) both contain a case of the spatial operator $-\chi_1(x)\Delta  + \chi_2(x)$ with 
\beq\label{E:per-gen} 
\chi_m(x+2\pi e_j)=\chi_m(x),  \quad \text{for all } x\in \R^d, j\in \{1,\dots,d\}, m\in \{1,2\}.
\eeq 
We also make the ellipticity assumption
\beq\label{E:ellip}
\chi_1(x)\geq \gamma >0, \ \chi_2(x) >0 \quad \text{for all }x\in \R^d.
\eeq

Applying the Bloch transformation produces 
$-\chi_1(x)|\nabla + \ri k|^2+ \chi_2(x)$. Thus, we consider the eigenvalue problem
\beq \label{E:Bloch-gen}
\cL(k)p_n(x,k)=\lambda_n(k)p_n(x,k), \ x\in \P, \ \cL(k):=-\chi_1(x)|\nabla + \ri k|^2 + \chi_2(x)
\eeq
with periodic boundary conditions. The operator $\mathcal{L}(k):H^2_\text{per}(\P)\rightarrow L^2_\text{per}(\P)$ is elliptic and
self adjoint in $L^2_{\chi_1}(\P,\C)$ equipped with the inner product $\langle f,g
\rangle_{L^2_{\chi_1}}=\int_\P
f(x)\overline{g}(x)\frac{1}{\chi_1(x)}\mathrm{d}x$. Because of
$\chi_1(x)\geq \gamma >0$ the induced norm
$\|\cdot\|_{L^2_{\chi_1}}$ is equivalent to the usual $L^2$ norm.
 The operator has a compact resolvent for each $k\in \B$. 
Hence, we conclude the existence of infinitely many real eigenvalues
$\lambda_n(k)$, $n\in\N$, with $\lambda_n(k)\rightarrow \infty$ for
$n\rightarrow\infty$. The spectrum $\text{spec}(-\chi_1 \Delta  + \chi_2)$ is real
and given by 
\beqq 
\text{spec}(-\chi_1 \Delta  + \chi_2)=\bigcup_{n\in\N, k\in\B}\lambda_n(k),
\eeqq 
see \cite[Chapter 3]{DLPSW2011}. We order
the eigenvalues by size $\lambda_n(k)\leq \lambda_{n+1}(k)$ for
$k\in\B$. The graph $(k,\lambda_{n}(k))_{n\in\N}$ is called the band
structure. As functions of $k$ the eigenvalues $\lambda_n(k)$ are
$1$-periodic in every component and analytic away from points of
higher multiplicity \cite[Sec. VII.2]{Kato_1995}. The eigenfunctions
$(p_n(\cdot,k))_{n\in \N}$ can be chosen to form an orthonormal
Schauder basis of $L_{\chi_1}^2(\P)$ for each $k\in \B$. As functions of $x$
they are $\P$-periodic and they are quasiperiodic in $k$:
\beq\label{E:quasiperpn} 
p_n(x,k+ e_j)=p_n(x,k)e^{-\ri x_j} 
\eeq 
for all $j\in\lbrace1,\cdots,d\rbrace$.

Our analysis requires certain regularity of the eigenfunctions. We provide here a simple $H^s$-regularity result based on a Fourier series analysis. Although the resulting statement may not be optimal, it allows fractional exponents $s>0$, which is suitable for us as we work with $\util(\cdot,k,t)$ in $\Hsper(\P)$.

\blem\label{L:reg}
Let $\chi_1,\chi_2 \in H^{s+d-2+\delta}_\text{per}(\P)$ with some $\delta >0$ satisfy \eqref{E:per-gen} and $\chi_1\geq \gamma>0$. Then $p_n(\cdot,k)\in \Hsper(\P)$.
\elem
\bpf
Because $\chi_1$ is bounded away from zero, we can rewrite \eqref{E:Bloch-gen} as
$$-|\nabla +\ri k|^2w = \psi(x) w, \ x \in \P, \text{ where }\psi(x):=\frac{1}{\chi_1(x)}(\lambda_n(k)-\chi_2(x))$$
for $w(\cdot):=p_n(\cdot,k)$. Expanding both $w$ and $\psi$ in the Fourier series
$$w(x)=\sum_{n\in \Z^d}W_n e^{\ri n\cdot x},  \ \psi(x) = \sum_{n\in \Z^d}\Psi_ne^{\ri n\cdot x},$$
the coefficient vector $\vec{W}:=(W_n)_{n\in \Z^d}$ satisfies
\beq\label{E:W}
(|n|^2+2k\cdot n+|k|^2)W_n=(\vec{\Psi}*\vec{W})_n,
\eeq
where $(\vec{\Psi}*\vec{W})_n=\sum_{j\in \Z^d}\Psi_{n-j}W_j$ and $\vec{\Psi}:=(\Psi_n)_{n\in \Z^d}$. 

For $w\in \Hsper(\P)$ we need to show $\vec{W}\in l^2_s(\Z^d)$, i.e. $\sum_{n\in \Z^d}|W_n|^2(1+|n|^{2s})<\infty.$ 
The idea is to first show $\vec{W}\in l^2_\tau(\Z^d), \tau\geq 0 \ \Rightarrow \ \vec{W}\in l^2_{\tau+2}(\Z^d)$ provided $\vec{\Psi}\in l^2_{d+\tau+\delta}(\Z^d)$ for  some $\delta >0$. Then, if $w$ solves \eqref{E:Bloch-gen}, then from \eqref{E:W} we get $(|n|^2W_n)_{n\in \Z^d}\in l^2_\tau(\Z^d)$, i.e. $\vec{W}\in l^2_{\tau+2}(\Z^d)$.

We have
$$
\begin{aligned}
\|\vec{\Psi}*\vec{W}\|_{l^2_\tau(\Z^d)}^2&=\sum_{n\in \Z^d}\big|\sum_{j\in\Z^d}\Psi_{n-j}W_j\big|^2(1+|n|^{2\tau}) \\
&\leq c\left(\sum_{n\in \Z^d}\big(\sum_{j\in\Z^d}|\Psi_{n-j}||W_j|(1+|n-j|^\tau)\big)^2+\sum_{n\in \Z^d}\big(\sum_{j\in\Z^d}|\Psi_{n-j}||W_j|(1+|j|^\tau)\big)^2\right)\\
&\leq c(\|\vec{\Psi}\|_{l^1_\tau(\Z^d)}^2\|\vec{W}\|_{l^2(\Z^d)}^2 + \|\vec{\Psi}\|_{l^1(\Z^d)}^2\|\vec{W}\|_{l^2_\tau(\Z^d)}^2),
\end{aligned}
$$
where the last step follows from Young's inequality for convolutions $\|f*g\|_{l^r}\leq \|f\|_{l^p}\|g\|_{l^q}$ for all $p,q,r\in [1,\infty)$ such that $1+\tfrac{1}{r}=\tfrac{1}{p}+\tfrac{1}{q}$. Next, we use the estimate 
$$|\Psi_n|\leq \frac{c}{|n|^a} \text{ for any } \psi\in H^a(\P), a>0,$$
which follows from $|\Psi_n||n|^a\leq \left(\sum_{n\in \Z^d}|\Psi_n|^2(1+|n|^{2a})\right)^{1/2}\leq c\|\psi\|_{H^a(\P)}.$ Because $\sum_{n\in \Z^d}|n|^{\tau-a}< \infty$ if and only if $a>d+\tau$, we get $\vec{\Psi}*\vec{W}\in l^2_\tau(\Z^d)$ if $\vec{W}\in l^2_{\tau}(\Z^d)$ and $\vec{\Psi}\in l^2_{d+\tau+\delta}(\Z^d)$ with $\delta>0$. The smoothness $w\in \Hsper(\P)$  can thus be concluded if $\chi_1,\chi_2 \in H_\text{per}^{d+s-2+\delta}(\P)$ with some $\delta>0$.
\epf

The asymptotic approximations below assume the This

As our estimates are performed in the Bloch variables and we often expand quantities in $k$, we use the Lipschitz continuity of the eigenfunction $p_{n_0}$:
\beq\label{E:Lipschitzpn}
\|p_{n_0}(\cdot,k)-p_{n_0}(\cdot,k_0)\|_{L^2(\P)}\leq L |k-k_0|,\quad L>0
\eeq
with some $L>0$ and for all $k$ in a neighborhood of $k_0$. For a simple eigenvalue $\lambda_{n_0}(k)$ at $k=k_0\in \B$ this is automatically satisfied, see \cite[Sec. VII.2]{Kato_1995}.

In Section \ref{S:esterror} we need also the Lipschitz continuity of $p_{n_0}$ with respect to $k$ in $H^n(\P), n>0$. In particular we use
\blem\label{L:Lip}
Let $q\in \N$, $\chi_1 \in H^{2q-2}_\text{per}(\P),\chi_2\in H^{2q-3}_\text{per}(\P), p_{n_0}(\cdot,k_0)\in H^{2q-1}_\text{per}(\P)$, and assume that $\lambda_{n_0}(k)$ is simple at $k=k_0$. Then
$$\left\|\cL^q(k)(p_{n_0}(\cdot,k)-p_{n_0}(\cdot,k_0))\right\|_{L^2(\P)}\leq L|k-k_0| \quad \text{for all } k\in \B,$$
where $L=L(\|p_{n_0}(\cdot,k_0)\|_{H^{2q-1}(\P)},q)$.
\elem
\bpf
We first note that
$$
\begin{aligned}
\cL(k)p_{n_0}(x,k_0) &= \cL(k_0)p_{n_0}(x,k_0)+\chi_1(k-k_0)^T(-2\ri \nabla p_{n_0}(x,k_0)+(k+k_0)p_{n_0}(x,k_0))\\
&=\lambda_{n_0}(k_0)p_{n_0}(x,k_0) + \alpha_1(x,k,k_0),
\end{aligned}
$$
where $\|\alpha_1(\cdot,k,k_0)\|_{L^2(\P)}\leq c|k-k_0|$.
Using a straightforward induction argument, one can show that under the conditions on $\chi_1,\chi_2$ and $ p_{n_0}(\cdot,k_0)$
$$\cL(k)^qp_{n_0}(x,k_0)=\lambda_{n_0}^q(k_0)p_{n_0}(x,k_0) + \alpha_q(x,k,k_0), \text{ with } \|\alpha_q(\cdot,k,k_0)\|_{L^2(\P)}\leq c|k-k_0|,$$
where $c=c(\|p_{n_0}(\cdot,k_0)\|_{H^{2q-1}(\P)})$. Then
$$
\begin{aligned}
\left\|\cL(k)^q(p_{n_0}(\cdot,k)-p_{n_0}(\cdot,k_0))\right\|_{L^2(\P)}&=\left\|\lambda_{n_0}^q(k)p_{n_0}(\cdot,k)-\lambda_{n_0}^q(k_0)p_{n_0}(\cdot,k_0)-\alpha_q(\cdot,k,k_0)\right\|_{L^2(\P)}\\
&\leq |\lambda_{n_0}^q(k)-\lambda_{n_0}^q(k_0)|\|p_{n_0}(\cdot,k)\|_{L^2(\P)} \\
& \ + |\lambda_{n_0}^q(k_0)|\|p_{n_0}(\cdot,k)-p_{n_0}(\cdot,k_0)\|_{L^2(\P)} + \|\alpha_q(\cdot,k,k_0)\|_{L^2(\P)}\\
&\leq L|k-k_0|
\end{aligned}
$$
with $L=L(\|p_{n_0}(\cdot,k_0)\|_{H^{2q-1}(\P)},q)$. In the last step we used \eqref{E:Lipschitzpn} and the Lipschitz continuity of $k\mapsto \lambda_{n_0}(k)$ \cite{CV1997}.
\epf


\subsection{Expansion in the Bloch eigenfunctions}

The completeness of the eigenfunctions allows to expand any function $\util(\cdot,k,t) \in H^s_\text{per}(\P)$. This in turn diagonalizes the operator $\cL(k)$, see equations \eqref{E:ODE-Un} and \eqref{E:ODE_Un2}.

We write
\beq\label{E:expansion}
\widetilde{u}(x,k,t)=\sum_{n\in\N}U_n(k,t)p_n(x,k), \quad U_n(k,t):=\langle \widetilde{u}(\cdot,k,t),p_n(\cdot,k)\rangle_{L^2_{\chi_1}(\P)}.
\eeq
Because of \eqref{E:perinxandk} and \eqref{E:quasiperpn}  we have
\beq\label{E:Un_per}
U_n(k+e_j,t)=U_n(k,t)  \text{ for all } 1\leq j \leq d, n\in\N, k \in \R^d.
\eeq
Next, we show that the diagonalization operator
\beq\label{E:Diag}
\mathcal{D}:L^1\left(\B,\Hsper(\P)\right)\rightarrow \cX(s):=L^1(\B,l^2_{s/d}),\quad \widetilde{u}\mapsto \vec{U}:=(U_n)_{n\in\N},\eeq
is an isomorphism.  This is shown in Lemma \ref{L:Diag}. Note that
$\|\vec{U}\|_{\cX(s)}:=\int_{\B}\|\vec{U}(k)\|_{l^2_{s/d}}\mathrm{d}k$
and
\beq \label{E:l2sd}
l^2_{s/d}:=\left\lbrace\vec{v}=(v_n)_{n\in\N}\in l^2(\mathbb{R}^d):\|\vec{v}\|^2_{l^2_{s/d}}=\sum_{n\in\N}n^{\frac{2s}{d}}|v_n|^2<\infty\right\rbrace.
\eeq
Hence, we can perform our estimates interchangeably in $L^1(\B,\Hsper(\P))$ and $\cX(s)$.

In order to show the isomorhism property of $\cD$ we first show the same for the $k-$dependent operator
\beq\label{E:ptDiag}
\mathcal{D}(k):\Hsper(\P)\rightarrow l^2_{s/d},\quad \widetilde{u}(\cdot,k)\mapsto \vec{U}(k):=(U_n(k))_{n\in\N}.
\eeq
\blem\label{L:ptDiag}
For all $s\geq 0$ and $k\in\B$ the operator $\mathcal{D}(k)$, defined in \eqref{E:ptDiag}, is an isomorphism between $\Hsper(\P)$ and $l^2_{s/d}$ and there exists a $\delta(k)>0$ with $\|\mathcal{D}(k)\|,\|\mathcal{D}^{-1}(k)\|\leq \delta(k)$. It is $\sup_{k\in\B}\delta(k)\leq C<\infty$.
\elem
\bpf
The proof is completely analogous to that for the case $d=1$ in Lemma 3.3 in \cite{BSTU06}. The main difference is that the asymptotics of the eigenvalues of $\cL(k)$ in $d$ dimensions are $C_1 n^{2/d}\leq \lambda_n(k) \leq C_2 n^{2/d}, \quad k \in \B, n \in \N$ with some $C_1,C_2>0$.
\epf
\blem\label{L:Diag}
For $s>\frac{d}{2}$ \beqq \mathcal{D}:L^1\left(\B,\Hsper(\P)\right) \rightarrow \cX(s), \quad \widetilde{u}\mapsto \vec{U}\eeqq is an isomorphism.
\elem
\bpf
We write $\vec{U}(k)=\mathcal{D}(k)\widetilde{u}(\cdot,k)$. Because of Lemma \ref{L:ptDiag} there are $c_1,c_2>0$ such that for all $k\in\B$ the estimates
\beqq \|\widetilde{u}(\cdot,k)\|_{H^s(\P)}\leq c_1\|\vec{U}(k)\|_{l^2_{s/d}}\eeqq
and
\beqq \|\vec{U}(k)\|_{l^2_{s/d}}\leq c_2\|\widetilde{u}(\cdot,k)\|_{H^s(\P)}\eeqq
are valid. It follows on the one hand
\beqq \|\widetilde{u}\|_{L^1(\B,\Hsper(\P))}\leq c_1\int_{\B}\|\vec{U}(k)\|_{l^2_{s/d}}\mathrm{d}k=c_1\|\vec{U}\|_{\cX(s)},\eeqq
on the other hand
\beqq \|\vec{U}\|_{\cX(s)}\leq c_2 \int_\B \|\widetilde{u}(\cdot,k)\|_{H^s(\P)}\mathrm{d}k=c_2\|\widetilde{u}\|_{L^1(\B,\Hsper(\P))}.\eeqq
\epf

Note that due to the periodicity of the coefficients $\vec{U}(\cdot,t)$
$$\|U_n(\cdot,t)\|_{L^1(\B)}=\|U_n(\cdot,t)\|_{L^1(\B+k_0)} \ \text{ for all } \ k_0\in \B$$
and thanks to the quasiperiodicity of the Bloch transform in $k$, we have
\beq \label{E:convol-shift}
\widetilde{u}*_\B\widetilde{v} = \widetilde{u}*_{\B+k_0}\widetilde{v}  \ \text{ for all } \ k_0\in \B
\eeq
because
$$
\begin{aligned}
(\widetilde{u}*_{\B+k_0}\widetilde{v})(x,k) &= \int_{\B+k_0}\widetilde{u}(x,k-l)\widetilde{v}(x,l)\dd l=\sum_{a\in \Z^d}\int_{(\B+k_0-a)\cap \B}\widetilde{u}(x,k-l-a)\widetilde{v}(x,l+a)\dd l\\
&= \sum_{a\in \Z^d}\int_{(\B+k_0-a)\cap \B}\widetilde{u}(x,k-l)\widetilde{v}(x,l)\dd l=(\widetilde{u}*_{\B}\widetilde{v})(x,k).
\end{aligned}
$$
We make use of these invariances of the norm and the convolution with respect to the $k_0-$shift in the proofs of the main results.


\section{The Gross-Pitaevskii Equation} \label{S:GP}

We consider first the Gross-Pitaevskii equation (GP)
\beq\label{E:GP} 
\ri \pa_t u+\Delta u -V(x) u -\sigma(x) |u|^2u=0,
\quad (x,T)\in \R^d\times \R 
\eeq 
with $d\in \N$ and the periodic
coefficients $V,\sigma$ such that 
\beq\label{E:per} V(x+2\pi
e_j)=V(x), \sigma(x+2\pi e_j)=\sigma(x) \quad \text{for all } x\in
\R^d, j\in \{1,\dots,d\}, \text{ and } V>0,\eeq 
where $e_j$ is the $j$-th Euclidean
unit vector in $\R^d$.

The formal asymptotic approximation is a wavepacket centered at a linear carrier wave, i.e. a Bloch wave. Bloch waves are determined from the Bloch eigenvalue problem \eqref{E:Bloch-gen}
with $\chi_1\equiv 1$ and $\chi_2:=V$. The eigenfunctions $(p_n(\cdot,k))_{n\in\N}$ are $L^2(\P)-$orthogonal in the Schr\"odinger case. The eigenvalues $\lambda_n$ play in the Schr\"odinger case the role of temporal frequencies, hence we use the notation
$$\omega_n(k):=\lambda_n(k), \ k\in \B$$ 
We select a $k_0\in \B$ and $n_0\in \N$, such that $\omega_{n_0}(k_0)$ is simple. The resulting Bloch wave is $p_{n_0}(x,k_0) e^{\ri(k_0\cdot x-\omega_0 t)}$, where $\omega_0:=\omega_{n_0}(k_0)$. This carrier propagates in the linear model ($\sigma \equiv 0$) at the group velocity
$$v_g:=\nabla \omega_{n_0}(k_0).$$
Hence, we expect the wavepacket to propagate at the same velocity and we make the following ansatz for an approximate solution of \eqref{E:GP}
\beq \label{E:uapp}
u_\text{app}(x,t)=\eps A(\eps(x-v_g t), \eps^2t) p_{n_0}(x,k_0) e^{\ri(k_0\cdot x-\omega_0 t)}.
\eeq
As shown in Sec. \ref{S:formal_der}, the effective equation for the envelope $A$ is
\beq\label{E:NLS}
\ri \pa_TA(X,T) +\frac{1}{2}\nabla \cdot (D^2\omega_{n_0}(k_0)\nabla A(X,T)) +\nu (|A|^2A)(X,T)=0, \quad (X,T) \in (\R^d\times \R),
\eeq
where $T:=\eps^2t, X:=\eps(x-v_g t), \nabla =\nabla_X$, and
$$\nu :=-\langle \sigma(\cdot)|p_{n_0}(\cdot,k_0)|^2p_{n_0}(\cdot,k_0),p_{n_0}(\cdot,k_0)\rangle_{L^2(\P)}.$$

Our main result is
\bthm\label{T:main}
Let $V\in H_\text{per}^{2d-2+\delta}(\P)$ and $\sigma \in H_\text{per}^{d+\delta}(\P)$ with some $\delta>0$ satisfy \eqref{E:per}. Assume that the Bloch eigenvalue $\omega_{n_0}:=\lambda_{n_0}$ of \eqref{E:Bloch-gen} with $\chi_1\equiv 1$ and $\chi_2:=V$ is simple at $k=k_0$ with the corresponding eigenfunction $p_{n_0}(\cdot,k_0)$. 
For every solution $A$ of \eqref{E:NLS} with the regularity $\widehat{A}\in C^1([0,T_0]; L^2(\R^d)\cap L^1_\beta(\R^d))$ for some $T_0>0$ and $\beta\geq 3$ with $\beta>2d$ there exist $\eps_0>0$ and $C>0$ such that for all $\eps\in (0,\eps_0)$ the solution $u$ of \eqref{E:GP} with initial data $u(x,0)=u_\text{app}(x,0)$ is continuous and satisfies
$$\sup_{x\in \R^d}|u(x,t)-u_\text{app}(x,t)|\leq C \eps^2 \quad \text{for all }t\in [0,T_0\eps^{-2}]$$
as well as the decay $u(x,t)\to 0$ for $|x|\to \infty$.
\ethm

The space $L^1_\beta(\R^d)$ is defined as
$$L^1_\beta(\R^d):=\{f\in L^1(\R^d):\int_{\R^d}|f(y)|(1+|y|^\beta) \dd y<\infty\}.$$

A large part of the analysis is carried out after applying $\cD\cT$ to the problem, i.e. we work with the Bloch expansion coefficients $U_n(k,t)$.
The expansion \eqref{E:expansion} leads to the infinite dimensional ODE-system
\beq \label{E:ODE-Un}
\ri\partial_t U_n(k,t)=\omega_n(k)U_n(k,t)+\langle\sigma(\cdot)(\widetilde{\bar{u}}*_{\B} \widetilde{u}*_{\B} \widetilde{u})(\cdot,k,t),p_n(\cdot,k)\rangle_{L^2(\P)}.\eeq

\subsection{Formal derivation of the effective NLS-equation}\label{S:formal_der}

The effective equation  \eqref{E:NLS} for the envelope comes as a necessary condition for making the residual small. We work in the Bloch variables, substitute the approximate ansatz \eqref{E:uapp} in \eqref{E:GP} and collect terms of equal $\eps$-power. In fact, \eqref{E:NLS} is derived carefully in Sec. \ref{S:res-calc} but we present here a shorter non-rigorous version, which gives more insight.

Using \eqref{E:xandtder}, \eqref{E:commper}, \eqref{E:conv}, and \eqref{E:quasiperpn}, the Bloch transformation of the approximate ansatz is
\beq \widetilde{u}_\text{app}(x,k,t)=\eps^{1-d}\sum_{m\in \Z^d}\widehat{A}\left(\frac{k-k_0+m}{\eps},T\right)p_{n_0}(x,k_0-m)e^{-\ri\alpha(k+m,k_0)t}, \quad (x,k,t)\in \P\times, \B\times\R\label{appAns}\eeq
where, for brevity, we use the shorthand notation
$$\alpha(k,k_0):=\omega_{n_0}(k_0)+(k-k_0)\cdot v_g.$$
Expanding in the orthonormal basis of the Bloch-functions $(p_n(\cdot,k))_{n\in\N}$
\beqq \widetilde{u}_\text{app}(x,k,t)=\sum_{n\in\N}U_n^\text{app}(k,t)p_n(x,k),\eeqq
we have the coefficients
\beq U_n^\text{app}(k,t)=\eps^{1-d}\sum_{m\in \Z^d}\widehat{A}\left(\frac{k-k_0+m}{\eps},T\right)\langle p_{n_0}(\cdot,k_0-m),p_n(\cdot,k)\rangle_{L^2(\P)}e^{-\ri\alpha(k+m,k_0)t}, \ (k,t)\in \B\times \R.\label{KoeffappAns}\eeq

Next, some non-rigorous approximations follow. These are all justified in Sec. \ref{S:res-calc} and \ref{S:res-est}.
First, due to the concentration of $\Ahat(\eps^{-1}(\cdot-k_0),\eps^2 t)$ near $k_0$ we set $\Ahat(\eps^{-1}(\cdot-k_0),\eps^2 t)$ to zero outside  $B_{\eps^{r}}(0)$ with some $r\in (0,1)$. Next, we approximate all $p_n(\cdot,k)$ by $p_n(\cdot,k_0)$. The compact support in $k$ implies that for $k\in \B+k_0$ we can reduce the sum over $m\in\Z^d$  to the one summand with $m=0$ if $\eps>0$ is sufficiently small. As a result
\begin{align}
\widetilde{u}_\text{app}(x,k,t)&\approx \eps^{1-d}\widehat{A}\left(\frac{k-k_0}{\eps},T\right)p_{n_0}(x,k_0)e^{-\ri\alpha(k,k_0)t}, \quad (x,k,t)\in \P\times (\B+k_0)\times\R, \label{E:util_appr}\\
U_n^\text{app}(k,t)&\approx\eps^{1-d}\widehat{A}\left(\frac{k-k_0}{\eps},T\right)\langle p_{n_0}(\cdot,k_0),p_n(\cdot,k_0)\rangle_{L^2(\P)}e^{-\ri\alpha(k,k_0)t}, \quad (k,t)\in (\B+k_0)\times\R,\\
&= \eps^{1-d}\widehat{A}\left(\frac{k-k_0}{\eps},T\right)\delta_{n,n_0}e^{-\ri\alpha(k,k_0)t}, \quad (k,t)\in (\B+k_0)\times\R. \label{E:Un_appr}
\end{align}
Outside $\B+k_0$ the approximation of $\widetilde{u}_\text{app}(x,k,t)$ in \eqref{E:util_appr} is defined via the quasiperiodicity in $k$, see \eqref{E:perinxandk}, and the approximation of $U_n^\text{app}(k,t)$ in \eqref{E:Un_appr} via the $1-$periodicity in $k$, see \eqref{E:Un_per}.

Next, we consider \eqref{E:ODE-Un}, where we Taylor-expand $\omega_{n_0}(k)$ near $k=k_0$ and on the compact support of the ansatz, i.e. in the $\eps^{r}$-neighborhood of $k_0$, we approximate
$$\omega_{n_0}(k)\approx \omega_{n_0}(k_0)+(k-k_0)^T\nabla\omega_{n_0}(k_0)+\frac{1}{2}(k-k_0)^TD^2\omega_{n_0}(k_0)(k-k_0).
$$

If we substitute these approximations in \eqref{E:ODE-Un}, all terms at $O(\eps^{1-d})$ and $O(\eps^{2-d})$ vanish due to $\omega_0= \omega_{n_0}(k_0)$ and $v_g = \nabla\omega_{n_0}(k_0)$. At $O(\eps^{3-d})$ we get on the left hand side
$$(\ri \pa_T\Ahat(\kappa,T) - \frac{1}{2}\kappa^T D^2\omega_{n_0}(k_0) \kappa \Ahat(\kappa,T))e^{-\ri \alpha t} -\eps^{d-3}\langle \sigma(\cdot)(\widetilde{\bar{u}}_\text{app}*_{\B} \widetilde{u}_\text{app}*_{\B} \widetilde{u}_\text{app})(\cdot,k,t),p_{n_0}(\cdot,k)\rangle_{L^2(\P)}$$
with 
$$\kappa :=\frac{k-k_0}{\eps}.$$
Approximating again $p_{n_0}(\cdot,k)$ by $p_{n_0}(\cdot,k_0)$, we get, analogously to \eqref{E:NL-fin-form}, for $k\in \B+k_0$
$$
\begin{aligned}
\langle \sigma(\cdot)&(\widetilde{\bar{u}}_\text{app}*_{\B} \widetilde{u}_\text{app}*_{\B} \widetilde{u}_\text{app})(\cdot,k,t),p_{n_0}(\cdot,k)\rangle_{L^2(\P)} \\
&\approx -\eps^{3-d}\nu e^{-\ri \alpha t}\int_{B_{2\eps^{r-1}}(0)}\int_{B_{\eps^{r-1}}(0)}\widehat{\overline{A}}(\kappa-h,T)\widehat{A}(h-l,T)\widehat{A}(l,T)\dd l \dd h\\
&=-\eps^{3-d}\nu e^{-\ri \alpha t}(\widehat{\overline{A}} *\widehat{A}*\widehat{A})(\kappa,T)
\end{aligned}
$$
This leads us to the condition
$$\ri \pa_T\Ahat - \frac{1}{2}\kappa^T D^2\omega_{n_0}(k_0) \kappa \Ahat + \nu \widehat{\overline{A}} *\widehat{A}*\widehat{A}=0, \quad (\kappa,T)\in (\R^n\times \R),$$
which is the effective NLS equation \eqref{E:NLS} in Fourier variables.


\subsection{Definition of the extended ansatz}

In order to produce a small enough residual, we modify the approximate ansatz $\widetilde{u}_\text{app}$ in (\ref{appAns}) to the, so called, extended ansatz $\widetilde{u}_\text{ext}$. After estimating the resulting residual and the approximation error $u_\text{ext}-u$, we also estimate $u_\text{app}-u_\text{ext}$ and use the triangle inequality to show the smallness of $u-u_\text{app}$. We use the following modifications of $\vec{U}^\text{app}$ in (\ref{KoeffappAns}). Firstly, we cut the $k$-support on $\B+k_0$ to a small neighborhood of $k_0$, which is motivated by the strong localization of $\widehat{A}(\eps^{-1}(\cdot-k_0),T)$ near $k_0$. Secondly, we replace in $U_{n_0}^\text{ext}$ the inner product $\langle p_{n_0}(\cdot,k_0-m),p_{n_0}(\cdot,k)\rangle_{L^2(\P)}$ by 1. On the support of $U_{n_0}^\text{ext}(\cdot,t)$ within $\B+k_0$ this incurs a small error. In $U_{n}^\text{ext}, n \neq n_0$ we denote the envelopes by $\eps^2\widehat{B}_n$.
Lastly, we add correction terms of higher $\eps$-order which are designed to cancel out those next-to-leading order (in $\eps$) terms in the residual terms which do not vanish after imposing the effective equation \eqref{E:NLS} on $A$.

For this we define
\beqq \widehat{A}_r\left(\kappa,T\right):=\cX_{B_{\eps^{r-1}}(0)}(\kappa)\widehat{A}\left(\kappa,T\right)\eeqq
with $r\in (0,1)$ to be chosen at a later stage of the proof. The extended ansatz is given by
\beq \widetilde{u}_\text{ext}(x,k,t)=\sum_{n\in\N}U_n^\text{ext}(k,t)p_n(x,k)\label{erwAns}\eeq
with
\beq
U_n^\text{ext}(k,t):=\eps^{1-d}\sum_{m\in\Z^d}\widehat{A}_n\left(\frac{k-k_0+m}{\eps},T\right)e^{-\ri\alpha(k+m,k_0)t}, \ (k,t)\in(\B+k_0)\times \R,\label{KoefferwAns}
\eeq
where
\beq \widehat{A}_{n_0}(\kappa,T):=\widehat{A}_r(\kappa,T),\quad \widehat{A}_{n\neq n_0}(\kappa,T):=\eps^2\widehat{B}_n(\kappa,T).\label{KoeffTr}\eeq
In order for $\widetilde{u}_\text{ext}$ to define the Bloch transform of a function, we need to extend the coefficients outside $\B+k_0$  periodically in $k$, i.e. $U_n^\text{ext}(k+e_j,t)=U_n^\text{ext}(k,t)$ for $j\in \{1, \dots d\}$ and all $n\in\N$. From the quasiperiodicity $p_n(x,k+e_j)=p_n(x,k)e^{-ix_j}$ we get for $\widetilde{u}_\text{ext}$ the $1-$quasiperiodicity in $k$ from \eqref{E:perinxandk}. The functions $\widehat{B}_n$ from (\ref{KoeffTr}) will be later chosen such that the residual of (\ref{erwAns}) in \eqref{E:ODE-Un} is sufficiently small. Due to the compact support of $\widehat{A}_r$ also the support of $\widehat{B}_n$ will be compact, namely
\beq\label{E:suppBn}
\text{supp}(\widehat{B}_n)\subset B_{3\eps^{r-1}}(0)
\eeq
as follows from \eqref{DefBn}.

The compact support of $\widehat{A}_r$ and $\widehat{B}_n$ implies that in (\ref{KoefferwAns}) for $k\in \B+k_0$ we can reduce the sum over $m\in\Z^d$  to the one summand with $m=0$ if $\eps>0$ is sufficiently small. Therefore, we have
\beq U_n^\text{ext}(k,t)=\eps^{1-d}\widehat{A}_n\left(\frac{k-k_0}{\eps},T\right)e^{-\ri\alpha(k,k_0)t}, \quad k\in\B+k_0, \ n \in \N.\label{KoefferwAnsneu}\eeq

\subsection{Calculation of the residual}\label{S:res-calc}

For the extended ansatz (\ref{erwAns}) we define the following residual corresponding to equation \eqref{E:ODE-Un}.
\beq  \label{E:resn}
Res_n(k,t):=-\ri\partial_t U_n^\text{ext}(k,t)+\omega_n(k)U_n^\text{ext}(k,t)+\langle\sigma(\cdot)(\widetilde{\widebar{u}}_\text{ext}*_{\B+k_0} \widetilde{u}_\text{ext}*_{\B+k_0} \widetilde{u}_\text{ext})(\cdot,k,t),p_n(\cdot,k)\rangle_{L^2(\P)}.
\eeq
Note that the convolutions in the nonlinear term can be written over $\B+k_0$ due to \eqref{E:convol-shift}.
The time derivative term is
\begin{align*}
-\ri\partial_t U_n^\text{ext}(k,t)=&-\eps^{1-d}\left[\omega_{n_0}(k_0)\widehat{A}_n\left(\kappa,T\right)e^{-\ri\alpha(k,k_0)t}\right]-\eps^{2-d}\left[\kappa^T\nabla_{k}\omega_{n_0}(k_0)\widehat{A}_n\left(\kappa,T\right)e^{-\ri\alpha(k,k_0)t}\right]\\
&- \eps^{3-d}\left[\ri\partial_T\widehat{A}_n\left(\kappa,T\right)e^{-\ri\alpha(k,k_0)t}\right]
\end{align*}
for all $n\in\N$, where $\kappa=\frac{k-k_0}{\eps}.$
The first two of these terms are clearly eliminated in $Res_{n_0}$ by rewriting $\omega_{n_0}(k)$  via its Taylor expansion. We write
\beq
\omega_{n_0}(k)=\omega_{n_0}(k_0)+(k-k_0)^Tv_g+\frac{1}{2}(k-k_0)^TD^2\omega_{n_0}(k_0)(k-k_0)+\phi(k),\label{Defphi}
\eeq
where $\phi(k)$ denotes the error in the Taylor approximation.

Next, we concentrate on the nonlinear term in \eqref{E:resn}. Because of $\widetilde{\widebar{u}}(x,k,t)=\widebar{\widetilde{u}}(x,-k,t)$ we have
\beqq \widetilde{\widebar{u}}_\text{ext}(x,k,t)=\sum_{n\in\N}\eps^{1-d}\widehat{\widebar{A}}_n\left(\frac{k+k_0}{\eps},T\right)e^{\ri\alpha(k,k_0)t}p_{n}(x,k), \quad k \in -\B-k_0,\eeqq
such that for $k\in \B+k_0$
\begin{equation} \label{nichtlinT}
\begin{aligned}
&\langle\sigma (\cdot)(\widetilde{\widebar{u}}_\text{ext}*_{\B+k_0}\widetilde{u}_\text{ext}*_{\B+k_0}\widetilde{u}_\text{ext})(\cdot,k,t),p_n(\cdot,k)\rangle_{L^2(\P)}\\
=& \eps^{3-3d}\sum_{a,b\in Z}\int_{\B+k_0}\int_{B_{\eps^r}(k_0)}\widehat{\widebar{A}}_{n_0}\left(\frac{k-h+k_0-a}{\eps},T\right)\widehat{A}_{n_0}\left(\frac{h-l-k_0+b}{\eps},T\right)\widehat{A}_{n_0}\left(\frac{l-k_0}{\eps},T\right)\\
& \qquad \qquad\times b_{n_0,n_0,n_0}^n(-k+h+a,h-l+b,l,k) \mathrm{d}l\mathrm{d}h e^{-\ri\alpha(k-a+b,k_0)t}+ R_{\eps^{5-d},n}(k,t),
\end{aligned}
\end{equation}
where
\begin{equation} \label{E:R-term}
\begin{aligned}
&R_{\eps^{5-d},n}(k,t):= \eps^{3-3d}\sum_{\substack{\lambda,\mu,\rho\in\N\\(\lambda,\mu,\rho)\neq(n_0,n_0,n_0)}}\sum_{a,b\in Z}\int_{\B+k_0}\int_{B_{3\eps^r}(k_0)}\widehat{\widebar{A}}_\lambda\left(\frac{k-h+k_0-a}{\eps},T\right)\\
&\times\widehat{A}_\mu\left(\frac{h-l-k_0+b}{\eps},T\right)\widehat{A}_\rho\left(\frac{l-k_0}{\eps},T\right)b_{\lambda,\mu,\rho}^{n}(-k+h+a,h-l+b,l,k) \mathrm{d}l\mathrm{d}h e^{-\ri\alpha(k-a+b,k_0)t},
\end{aligned}
\end{equation}
with
\beqq
Z:=\left\lbrace z\in\Z^d : |z_j|\leq 1, 1\leq j\leq d\right\rbrace
\label{DefZ}
\eeqq
and
\beqq
b_{\lambda,\mu,\rho}^n(c,d,e,f):=\langle \sigma(\cdot)\overline{p_\lambda}(\cdot,c)p_\mu(\cdot,d)p_\rho(\cdot,e),p_n(\cdot,f)\rangle_{L^2(\P)}, \ n,\lambda,\mu,\rho \in \N.
\eeqq
Note that the domain for the $l-$integrals has been reduced from $\B+k_0$ to $B_{\eps^r}(k_0)$ (resp. $B_{3\eps^r}(k_0)$)  due to the compact support of $\widehat{A}_{n}$. Note that in \eqref{nichtlinT} and in \eqref{E:R-term} we have used $p_n(x,k)=\overline{p_n}(x,-k)$.

The term $R_{\eps^{5-d},n}$ collects higher order nonlinearity terms, i.e. those where at least one of $\lambda,\mu,\rho$ is not $n_0$, such that all terms in $R_{\eps^{5-d},n}$  are of the formal order $\eps^{5-d}$. This power of $\eps$ results after the variable changes $\tilde{l}:=(l-k_0)/\eps$ and $\tilde{h}:=(h-2k_0)/\eps$ in the integral and because at least one of the $\widehat{A}_i$'s is $\eps^2\widehat{B}_i$. Note that the summation over $a,b\in Z$ is necessary as $h-l$, resp. $k-h$ do not always lie in $\B+k_0$ resp. $\B-k_0$  if $k,h\in \B+k_0, l\in B_{\eps^r}(k_0)$. For $k\in \B+k_0$ it is, however, $a=b$ due to the support of $\widehat{A}_n$. This can be seen as follows. In \eqref{nichtlinT} due to $\text{supp}(\widehat{A}_{n_0})\subset B_{\eps^{r-1}}(0), l\in B_{\eps^r}(k_0)$, and $h\in \B+k_0$ it is $h-l-k_0+b \in  B_{\eps^r}(0)$ if and only if $h\in B_{2\eps^r}(2k_0-b)$. For the argument of $\widehat{\widebar{A}}_{n_0}$ we have due to $h\in B_{2\eps^r}(2k_0-b)$ analogously $k\in B_{3\eps^r}(k_0+a-b)$. Hence, for $k\in \B+k_0$, it must be $a=b$ because otherwise $|a-b|\geq 1$ such that $k_0+a-b\notin \overline{\B}+k_0$. Therefore, the double integral in \eqref{nichtlinT} is
\begin{align*}
\sum_{a \in Z} \int_{B_{2\eps^r}(2k_0-a)\cap(\B+k_0)}\int_{B_{\eps^r}(k_0)}&\widehat{\widebar{A}}_{n_0}\left(\frac{k-h+k_0-a}{\eps},T\right)\widehat{A}_{n_0}\left(\frac{h-l-k_0+a}{\eps},T\right)\widehat{A}_{n_0}\left(\frac{l-k_0}{\eps},T\right)\notag\\
&\times b_{n_0,n_0,n_0}^{n}(-k+h+a,h-l+a,l,k) \mathrm{d}l\mathrm{d}h e^{-\ri\alpha(k,k_0)t}.
\end{align*}
After the substitution $h':=h+a$ it is easy to see that the sum over $a$ simply produces the $h'-$integral over the full $B_{2\eps^r}(2k_0)$. The calculation for $R_{\eps^{5-d},n}$ is completely analogous. Renaming $h'$ to $h$ again, we get for $k\in \B+k_0$
\begin{align}
\langle\sigma (\cdot)(\widetilde{\widebar{u}}_\text{ext}*_{\B+k_0}&\widetilde{u}_\text{ext}*_{\B+k_0}\widetilde{u}_\text{ext})(\cdot,k,t),p_n(\cdot,k)\rangle_{L^2(\P)}\notag\\
=& \eps^{3-3d}\int_{B_{2\eps^r}(2k_0)}\int_{B_{\eps^r}(k_0)}\widehat{\widebar{A}}_{n_0}\left(\frac{k-h+k_0}{\eps},T\right)\widehat{A}_{n_0}\left(\frac{h-l-k_0}{\eps},T\right)\widehat{A}_{n_0}\left(\frac{l-k_0}{\eps},T\right)\notag\\
& \qquad \qquad\times b_{n_0,n_0,n_0}^n(-k+h,h-l,l,k) \mathrm{d}l\mathrm{d}h e^{-\ri\alpha(k,k_0)t}+ R_{\eps^{5-d},n}(k,t),
\label{nichtlinT2}
\end{align}
where
\begin{align}
R_{\eps^{5-d},n}(k,t)=& \eps^{3-3d}\int_{B_{6\eps^r}(2k_0)}\int_{B_{3\eps^r}(k_0)}\sum_{\substack{\lambda,\mu,\rho\in\N\\(\lambda,\mu,\rho)\neq(n_0,n_0,n_0)}}\hspace{-0.4cm}\widehat{\widebar{A}}_\lambda\left(\frac{k-h+k_0}{\eps},T\right)\widehat{A}_\mu\left(\frac{h-l-k_0}{\eps},T\right)\notag\\
&\qquad\qquad \times\widehat{A}_\rho\left(\frac{l-k_0}{\eps},T\right)b_{\lambda,\mu,\rho}^{n}(-k+h,h-l,l,k) \mathrm{d}l\mathrm{d}h~e^{-\ri\alpha(k,k_0)t}.\label{R-term}
\end{align}

We further rewrite the leading order nonlinear term using the variables
\beqq \widetilde{l}:=\frac{l-k_0}{\eps}, \quad \widetilde{h}:=\frac{h-2k_0}{\eps}.\eeqq
Clearly  $\widetilde{l}\in B_{\eps^{r-1}}(0)$ and $\widetilde{h} \in B_{2\eps^{r-1}}(0)$ and we have for $k\in \B+k_0$
\beq \label{E:NL-fin-form}
\begin{aligned}
\langle\sigma (\cdot)(\widetilde{\widebar{u}}_\text{ext}*_{\B+k_0}\widetilde{u}_\text{ext}*_{\B+k_0}\widetilde{u}_\text{ext})(\cdot,k,t),&p_n(\cdot,k)\rangle_{L^2(\P)}\\
= \eps^{3-d}\int_{B_{2\eps^{r-1}}(0)}\int_{B_{\eps^{r-1}}(0)}&\widehat{\widebar{A}}_{n_0}\left(\kappa-\widetilde{h},T\right)\widehat{A}_{n_0}\left(\widetilde{h}-\widetilde{l},T\right)\widehat{A}_{n_0}\left(\widetilde{l},T\right)\\
&\times \widetilde{b}_{n_0,n_0,n_0}^n(\kappa,\widetilde{h},\widetilde{l}) \mathrm{d}\widetilde{l}\mathrm{d}\widetilde{h}e^{-\ri\alpha(k,k_0)t}
+ R_{\eps^{5-d},n}(k,t),
\end{aligned}
\eeq
where
\beqq \widetilde{b}_{\lambda,\mu,\rho}^n(\kappa,\widetilde{h},\widetilde{l}):=b_{\lambda,\mu,\rho}^n(\eps(\widetilde{h}-\kappa)+k_0,\eps(\widetilde{h}-\widetilde{l})+k_0,\eps\widetilde{l}+k_0,\eps\kappa +k_0).\eeqq

Next, in order to recover the effective NLS equation \eqref{E:NLS}, we approximate in the case $n=n_0$ the function $\widetilde{b}^{n_0}_{n_0,n_0,n_0}(\kappa,\widetilde{h},\widetilde{l})$ by its value at $\eps=0$ and the convolution of $\widehat{A}_{n_0}$ by the full space convolution of $\widehat{A}$. For that we recall from \eqref{E:NLS}
$$\nu=- b_{n_0,n_0,n_0}^{n_0}(k_0,k_0,k_0,k_0)$$
and we write
\begin{align*}
\langle\sigma (\cdot)(\widetilde{\widebar{u}}_\text{ext}&*_{\B+k_0}\widetilde{u}_\text{ext}*_{\B+k_0}\widetilde{u}_\text{ext})(\cdot,k,t),p_{n_0}(\cdot,k)\rangle_{L^2(\P)}\notag\\
= & -\nu\eps^{3-d}\chi_{B_{\eps^{r-1}}(0)}(\kappa)\left(\widehat{\widebar{A}}*_{\R^d}\widehat{A}*_{\R^d}\widehat{A}\right)\left(\kappa,T\right)e^{-\ri\alpha(k,k_0)t} + \psi(k,t)+R_{\eps^{5-d},n}(k,t),
\end{align*}
 where
\beq\label{Defpsin}
\begin{aligned}
\psi(k,t):=& -\nu\eps^{3-d}\left(\widehat{\widebar{A}}_{n_0}*_{B_{2\eps^{r-1}}(0)}\widehat{A}_{n_0}*_{B_{\eps^{r-1}}(0)}\widehat{A}_{n_0}-\chi_{B_{\eps^{r-1}}(0)}\widehat{\widebar{A}}*_{\R^d}\widehat{A}*_{\R^d}\widehat{A}\right)\left(\kappa,T\right)e^{-\ri\alpha(k,k_0)t}\\
& +\eps^{3-d}\int_{B_{2\eps^{r-1}}(0)}\int_{B_{\eps^{r-1}}(0)}\widehat{\widebar{A}}_{n_0}\left(\kappa-\widetilde{h},T\right)\widehat{A}_{n_0}\left(\widetilde{h}-\widetilde{l},T\right)\widehat{A}_{n_0}\left(\widetilde{l},T\right)\\
&\qquad\qquad\qquad\qquad \times \left(\widetilde{b}_{n_0,n_0,n_0}^{n_0}(\kappa,\widetilde{h},\widetilde{l})-b_{n_0,n_0,n_0}^{n_0}(k_0,k_0,k_0,k_0)\right)\mathrm{d}\widetilde{l}\mathrm{d}\widetilde{h}~e^{-\ri\alpha(k,k_0)t}.
\end{aligned}
\eeq

Collecting now all terms of the residual (\ref{E:resn}) and using for $n=n_0$ the effective NLS equation \eqref{E:NLS}, we obtain
\begin{align}
Res_{n_0}(k,t)=&\eps^{3-d}\left[-\ri\partial_T\widehat{A}_r\left(\kappa,T\right)+\tfrac{1}{2} \kappa^T D^2\omega_{n_0}(k_0)\kappa \widehat{A}_r\left(\kappa,T\right)\right. \notag\\
&\hspace{-2cm}\left.- \nu\chi_{B_{\eps^{r-1}}(0)}(\kappa)\left(\widehat{\widebar{A}}*_{\R^d}\widehat{A}*_{\R^d}\widehat{A}\right)\left(\kappa,T\right)\right]e^{-\ri\alpha(k,k_0)t}+\psi(k,t)+\phi(k)U_{n_0}^\text{ext}(k,t)+R_{\eps^{5-d},n_0}(k,t)\notag\\
=& \psi(k,t)+\phi(k)U_{n_0}^\text{ext}(k,t)+R_{\eps^{5-d},n_0}(k,t)\label{Resn0}
\end{align}
and
\begin{align*}
Res_{n\neq n_0}(k,t)=& \eps^{3-d} (-\omega_{n_0}(k_0)+\omega_n(k))\widehat{B}_n\left(\kappa,T\right)e^{-\ri\alpha(k,k_0)t}\\
-& \eps^{3-d} (k-k_0)^T \nabla_{k}\omega_{n_0}(k_0)\widehat{B}_n\left(\kappa,T\right)e^{-\ri\alpha(k,k_0)t} \\
+&\eps^{3-d}\int_{B_{2\eps^{r-1}}(0)}\int_{B_{\eps^{r-1}}(0)}\widehat{\widebar{A}}_{n_0}\left(\kappa-\widetilde{h},T\right)\widehat{A}_{n_0}\left(\widetilde{h}-\widetilde{l},T\right)\widehat{A}_{n_0}\left(\widetilde{l},T\right)\\
&\hspace{2cm}\times \widetilde{b}_{n_0,n_0,n_0}^n(\kappa,\widetilde{h},\widetilde{l}) \mathrm{d}\widetilde{l}\mathrm{d}\widetilde{h}e^{-\ri\alpha(k,k_0)t}\\
+&R_{\eps^{5-d},n}(k,t)-\eps^{5-d}\ri\partial_T \widehat{B}_n\left(\frac{k-k_0}{\eps},T\right)e^{-\ri\alpha(k,k_0)t}.
\end{align*}

We complete the definition of the extended ansatz in (\ref{erwAns}), (\ref{KoefferwAns}), and \eqref{KoeffTr} by the subsequent choice of
\beq\label{DefBn}
\begin{aligned}
\widehat{B}_n\left(\frac{k-k_0}{\eps},T\right):= \frac{1}{\eta_n(k,k_0)}\int_{B_{2\eps^{r-1}}(0)}\int_{B_{\eps^{r-1}}(0)}&\widehat{\widebar{A}}_r\left(\kappa-\widetilde{h},T\right)\widehat{A}_r\left(\widetilde{h}-\widetilde{l},T\right)\\
&\times \widehat{A}_r\left(\widetilde{l},T\right)\widetilde{b}_{n_0,n_0,n_0}^n(\kappa,\widetilde{h},\widetilde{l}) \mathrm{d}\widetilde{l}\mathrm{d}\widetilde{h}
\end{aligned}
\eeq
with
\beq \eta_n(k,k_0):=\omega_{n_0}(k_0)+(k-k_0)^Tv_g-\omega_n(k).\label{Defetan}\eeq
This choice eliminates all terms of formal order $\eps^{3-d}$ in $Res_{n\neq n_0}$ such that
\beq Res_{n\neq n_0}(k,t)=R_{\eps^{5-d},n}(k,t)-\eps^{5-d}\ri\partial_T \widehat{B}_n\left(\frac{k-k_0}{\eps},T\right)e^{-\ri\alpha(k,k_0)t}.\label{Resn}\eeq

The above definition of $\widehat{B}_n$ satisfies the condition $\text{supp}(\widehat{B}_n) \subset B_{3\eps^{r-1}}(0)$ in \eqref{E:suppBn} due to the double convolution structure. Also note that due to our assumption of the simpleness of $\omega_{n_0}$ at $k=k_0$, the denominator $\eta_n(k,k_0)$ is bounded away from zero for all $k\in B_{3\eps^{r}}(k_0)$, $n\neq n_0$ and $\eps>0$ small enough.

\subsection{Estimation of the residual}\label{S:res-est}

We estimate $(Res_n(k,t))_{n\in\N}$ in $\cX(s)=L^1(\B,l^2_{s/d})$. For that we first show the smallness (as $\eps\to 0$) of the leading order parts $\|\psi(\cdot,t)\|_{L^1(\B)}$ and $\|\phi(\cdot)U_{n_0}^\text{ext}(\cdot,t)\|_{L^1(\B)}$. Next, for $\eps^{5-d}\partial_T \widehat{B}_n\left(\tfrac{\cdot-k_0}{\eps},T\right)$ and $R_{\eps^{5-d},n}(\cdot,t)$ we first estimate each of these in $L^1(\B)$ by a constant $c_n(\eps)$ such that $(c_n(\eps))_n \in l^2_{s/d}$. Note that $\eps$-independent constants are denoted by $C$ and their meaning often changes from one line to the next.

We have
\beq\label{Abschpsi}
\begin{aligned}
&\|\psi(\cdot,t)\|_{L^1(\B+k_0)}=\eps^d\|\psi(\eps\cdot+k_0,t)\|_{L^1\left(\frac{\B}{\eps}\right)}\\
\leq & |\nu|\eps^{3}\left\|\left(\widehat{\widebar{A}}_r*_{B_{2\eps^{r}}(0)}\widehat{A}_r*_{B_{\eps^r}(0)}\widehat{A}_r-\chi_{B_{\eps^{r-1}}(0)}\widehat{\widebar{A}}*_{\R^d}\widehat{A}*_{\R^d}\widehat{A}\right)\left(\cdot,T\right)\right\|_{L^1\left(\frac{\B}{\eps}\right)}\\
& +\eps^{3}\int_{B_{3\eps^{r-1}}(0)}\int_{B_{2\eps^{r-1}}(0)}\int_{B_{\eps^{r-1}}(0)}\left|\widehat{\widebar{A}}_r\left(\kappa-\widetilde{h},T\right)\widehat{A}_r\left(\widetilde{h}-\widetilde{l},T\right)\widehat{A}_r\left(\widetilde{l},T\right)\right|\\
& \qquad\qquad\qquad\qquad\qquad\qquad\times \left|\widetilde{b}_{n_0,n_0,n_0}^{n_0}(\kappa,\widetilde{h},\widetilde{l})-b_{n_0,n_0,n_0}^{n_0}(k_0,k_0,k_0,k_0)\right|\mathrm{d}\widetilde{l}\mathrm{d}\widetilde{h}\mathrm{d}\kappa,
\end{aligned}
\eeq
in which the multilinearity of $b_{n_0,n_0,n_0}^{n_0}$ provides
\begin{align*}
& |b_{n_0,n_0,n_0}^{n_0}(k_0+\eps(\widetilde{h}-\kappa),\eps(\widetilde{h}-\widetilde{l})+k_0,\eps\widetilde{l}+k_0,\eps\kappa+k_0)-b_{n_0,n_0,n_0}^{n_0}(k_0,k_0,k_0,k_0)|\\
\leq &|b_{n_0,n_0, n_0}^{n_0}(k_0+\eps(\widetilde{h}-\kappa),k_0+\eps(\widetilde{h}-\widetilde{l}),k_0+\eps\widetilde{l},k_0+\eps\kappa) \\
&\hspace{3cm}-b_{n_0,n_0, n_0}^{n_0}(k_0,k_0+\eps(\widetilde{h}-\widetilde{l}),k_0+\eps\widetilde{l},k_0+\eps\kappa)|\\
&+\cdots +|b_{n_0,n_0, n_0}^{n_0}(k_0,k_0,k_0,k_0+\eps\kappa)-b_{n_0,n_0, n_0}^{n_0}(k_0,k_0,k_0,k_0)|.
\end{align*}
In every summand we use the Cauchy-Schwarz-inequality, the Lipschitz-continuity of the Bloch-functions \eqref{E:Lipschitzpn} and the algebra-property of the $H^s$-norm for $s>d/2$. For example, provided $\sigma,p_{n_0}(\cdot,k_0)\in \Hsper(\P)$, the last summand yields
\begin{align*}
&|b_{n_0,n_0, n_0}^{n_0}(k_0,k_0,k_0,k_0+\eps\kappa)-b_{n_0,n_0, n_0}^{n_0}(k_0,k_0,k_0,k_0)|\\
\leq & C \eps |\kappa | L \|\sigma(\cdot)p_{n_0}(\cdot,k_0)p_{n_0}(\cdot,k_0)p_{n_0}(\cdot,k_0)\|_{H^{s}(\P)}\leq  C \eps |\kappa| L \|\sigma(\cdot)\|_{H^{s}(\P)}\|p_{n_0}(\cdot,k_0)\|^3_{H^{s}(\P)}\leq C \eps  |\kappa|.
\end{align*}
The $H^{s}$-regularity of $p_{n_0}(\cdot,k_0)$ follows if $V\in H^{a}_\text{per}(\P)$ with $a>s+d-2$, see Lemma \ref{L:reg}.
Altogether we get
\beqq \left|\widetilde{b}_{n_0,n_0,n_0}^{n_0}(\kappa,\widetilde{h},\widetilde{l})-b_{n_0,n_0,n_0}^{n_0}(k_0,k_0,k_0,k_0)\right| \leq C\eps \left(|\kappa - \widetilde{h}|+|\widetilde{h}-\widetilde{l}|+|\widetilde{l}|+|\kappa|\right).\eeqq
Hence, the second term in (\ref{Abschpsi}) is bounded by
\begin{align*}
&C\eps^4 \int_{B_{3\eps^{r-1}}(0)} \bigg(h(\kappa)\left|\left(\widehat{\widebar{A}}_r*_{B_{2\eps^{r-1}}(0)}\widehat{A}_r*_{B_{\eps^{r-1}}(0)}\widehat{A}_r\right)\left(\kappa,T\right)\right| + \left|\left(h\widehat{\widebar{A}}_r*_{B_{2\eps^{r-1}}(0)}\widehat{A}_r*_{B_{\eps^{r-1}}(0)}\widehat{A}_r\right)\left(\kappa,T\right)\right|\\
&\qquad\qquad+ \left|\left(\widehat{\widebar{A}}_r*_{B_{2\eps^{r-1}}(0)}h\widehat{A}_r*_{B_{\eps^{r-1}}(0)}\widehat{A}_r\right)\left(\kappa,T\right)\right|+ \left|\left(\widehat{\widebar{A}}_r*_{B_{2\eps^{r-1}}(0)}\widehat{A}_r*_{B_{\eps^{r-1}}(0)}h\widehat{A}_r\right)\left(\kappa,T\right)\right|\bigg)\mathrm{d}\kappa,
\end{align*}
where $h(\kappa):=|\kappa|$. This can be estimated with Young's inequality for convolutions by
\beq\label{E:psi1} 
C\eps^{4}\|\widehat{A}_r(\cdot,T)\|_{L_1^1(B_{\eps^{r-1}}(0))}\|\widehat{A}_r(\cdot,T)\|^2_{L^1(B_{\eps^{r-1}}(0))} \leq C\eps^{4} \|\widehat{A}(\cdot,T)\|^3_{L^1_1(\R^d)}.
\eeq
For an estimate of the first term  in (\ref{Abschpsi}) we denote by $-\widetilde{A}$ the tail of $\widehat{A}$, i.e.
\beq \widetilde{A}(\kappa,T):=\left(\cX_{B_{\eps^{r-1}}(0)}(\kappa)-1\right)\widehat{A}(\kappa,T).\label{ASchlange}\eeq
Obviously,
$$\widehat{A}_r=\widehat{A}+\widetilde{A}.$$

The smallness of the tail is demonstrated by the following calculation.
\begin{align}
&\|\widetilde{A}(\cdot,T)\|_{L^1(\R^d)}=\|(\chi_{B_{\eps^{r-1}}(0)}(\cdot)-1)\widehat{A}(\cdot,T)\|_{L^1(\R^d)} \notag\\
&=\|(\chi_{B_{\eps^{r-1}}(0)}(\cdot)-1)(1+|\cdot|)^{-\beta}(1+|\cdot|)^\beta\widehat{A}(\cdot,T)\|_{L^1(\R^d)}\notag\\
&\leq  \|(1+|\cdot|)^\beta\widehat{A}(\cdot,T)\|_{L^1(\R^d)} \sup_{\kappa \in \R^d}|(\chi_{B_{\eps^{r-1}}(0)}(\kappa)-1)(1+|\kappa|)^{-\beta}| \notag\\
&=  \|\widehat{A}(\cdot,T)\|_{L_\beta^1(\R^d)} \sup_{\kappa \in \R^d,|\kappa|>\eps^{r-1}}|(1+|\kappa|)^{-\beta}| \leq \eps^{(1-r)\beta}\|\widehat{A}(\cdot,T)\|_{L_\beta^1(\R^d)}
\label{L1beta}
\end{align}
if $\beta>0$ and $\widehat{A}(\cdot,T)\in L^1_\beta(\R^d)$.

We esimate the first term in  (\ref{Abschpsi}) via
$$ \left\|\left(\widehat{\widebar{A}}_r*_{B_{2\eps^{r-1}}(0)}\widehat{A}_r*_{B_{\eps^{r-1}}(0)}\widehat{A}_r\right)\left(\cdot,T\right)-\chi_{B_{\eps^{r-1}}(0)}(\cdot)\left(\widehat{\widebar{A}}*_{\R^d}\widehat{A}*_{\R^d}\widehat{A}\right)\left(\cdot,T\right)\right\|_{L^1\left(\frac{\B}{\eps}\right)} \leq I_1+I_2,$$
where
$$
\begin{aligned}
I_1&:=\left\|\left(\widehat{\widebar{A}}_r*_{B_{2\eps^{r-1}}(0)}\widehat{A}_r*_{B_{\eps^{r-1}}(0)}\widehat{A}_r-\widehat{\widebar{A}}*_{\R^d}\widehat{A}*_{\R^d}\widehat{A}\right)\left(\cdot,T\right)\right|\|_{L^1\left(\frac{\B}{\eps}\right)}\\
I_2&:=\left\|\chi_{\eps^{-1}\B\setminus B_{\eps^{r-1}}(0)}(\cdot)\left(\widehat{\widebar{A}}*_{\R^d}\widehat{A}*_{\R^d}\widehat{A}\right)\left(\cdot,T\right)\right\|_{L^1\left(\frac{\B}{\eps}\right)}.
\end{aligned}
$$
For $I_1$ we have
\begin{align}
I_1= & \left\|\left((\widehat{\widebar{A}}+\widetilde{\widebar{A}})*_{\R^d}(\widehat{A}+\widetilde{A})*_{\R^d}(\widehat{A}+\widetilde{A})-\widehat{\widebar{A}}*_{\R^d}\widehat{A}*_{\R^d}\widehat{A}\right)\left(\cdot,T\right)\right\|_{L^1\left(\frac{\B}{\eps}\right)}\notag\\
\leq & \left\|\left(\widetilde{\widebar{A}}*_{\R^d} \widehat{A}*_{\R^d} \widehat{A}\right)\left(\cdot,T\right)\right\|_{L^1(\R^d)} +\left\|\left(\widehat{\widebar{A}}*_{\R^d} \widetilde{A}*_{\R^d} \widehat{A}\right)\left(\cdot,T\right)\right\|_{L^1(\R^d)}\notag\\
&+\left\|\left(\widehat{\widebar{A}}*_{\R^d} \widehat{A}*_{\R^d} \widetilde{A}\right)\left(\cdot,T\right)\right\|_{L^1(\R^d)}+ \text{ terms  quadratic or cubic in $\widetilde{A}$}\notag\\
\leq & C\left(\|\widetilde{A}(\cdot,T)\|_{L^1(\R^d)}\|\widehat{A}(\cdot,T)\|^2_{L^1(\R^d)}+\|\widetilde{A}(\cdot,T)\|^2_{L^1(\R^d)}\|\widehat{A}(\cdot,T)\|_{L^1(\R^d)}+\|\widetilde{A}(\cdot,T)\|^3_{L^1(\R^d)}\right)\notag\\
\leq & C \eps^{(1-r)\beta}\|\widehat{A}(\cdot,T)\|^3_{L_\beta^1(\R^d)}.\label{E:psi2}
\end{align}
For $I_2$ we first define $\widetilde{A}^{(1/3)}(\kappa,T):=\left(1-\chi_{B_{\frac{\eps^{r-1}}{3}}(0)}(\kappa)\right)\widehat{A}(\kappa,T).$ Due to \eqref{L1beta} we have
$$\|\widetilde{A}^{(1/3)}(\cdot,T)\|_{L^1(\R^d)}\leq \left(\tfrac{1}{3}\right)^{(1-r)\beta}\eps^{(1-r)\beta}\|\widehat{A}(\cdot,T)\|_{L_\beta^1(\R^d)}.$$
Next, because
$$\left((\chi_{B_{\eps^{r-1}/3}(0)}\widehat{\overline{A}})*_{\R^d}(\chi_{B_{\eps^{r-1}/3}(0)}\widehat{A})*_{\R^d}(\chi_{B_{\eps^{r-1}/3}(0)}\widehat{A})\right)(\kappa,T)=0 \text{  for } \kappa \in \eps^{-1}\B\setminus B_{\eps^{r-1}(0)}, $$
we have
\begin{align}
I_2&=\left\|\chi_{\eps^{-1}\B\setminus B_{\eps^{r-1}}(0)}(\cdot)\left(\widetilde{\widebar{A}}^{(1/3)}*_{\R^d}\widehat{A}*_{\R^d}\widehat{A}\right.\right.\notag\\
&\qquad +\left.\left.\chi_{B_{\eps^{r-1}/3}(0)}\widehat{\overline{A}}*_{\R^d}\widehat{A}*_{\R^d}(\widetilde{A}^{(1/3)}+\chi_{B_{\eps^{r-1}/3}(0)}\widehat{A})\right)\left(\cdot,T\right)\right\|_{L^1\left(\frac{\B}{\eps}\right)}\notag\\
&\leq C\eps^{(1-r)\beta}\|\widehat{A}(\cdot,T)\|^3_{L^1_\beta(\R^d)}.\label{E:psi3}
\end{align}

In summary using \eqref{E:psi1},\eqref{E:psi2}, and \eqref{E:psi3},
\beqq
\|\psi(\cdot,t)\|_{L^1(\B)} \leq C\eps^4\left(\|\widehat{A}(\cdot,T)\|^3_{L^1_1(\R^d)}+\|\widehat{A}(\cdot,T)\|^3_{L_\beta^1(\R^d)}\eps^{(1-r)\beta-1}\right).
\eeqq
For an optimal estimate we need to set the free parameter $\beta$ so that $(1-r)\beta-1\geq 0$, i.e.
\beq\label{E:beta-r-rel}
\beta\geq \frac{1}{1-r}.
\eeq
Because $r\in (0,1)$, it is then $\beta >1$ and we get
\beq\label{E:psi-est}
\|\psi(\cdot,t)\|_{L^1(\B)} \leq C\eps^4\|\widehat{A}(\cdot,T)\|^3_{L^1_\beta(\R^d)}
\eeq
under the regularity condition $\widehat{A}(\cdot,T) \in L^1_\beta(\R^d)$.

For an estimate of $\|\phi(\cdot)U_{n_0}^\text{ext}(\cdot,t)\|_{L^1(\B)}$ we first prove the following lemma
\blem\label{L:gD-est}
 Let $a,b\geq 0$, $D\in L^1_{a+b}(\R^d)$, and $g\in L^1(\R^d)$ with $|g(k)|\leq C |k-k_0|^b$ for all $k\in \R^d$ and some $C>0$. Then
\beqq \|g(\eps\cdot +k_0)D(\cdot)\|_{L_a^1(\R^d)}\leq C\eps^b\|D\|_{L_{a+b}^1(\R^d)}.\eeqq
\elem
\bpf
The proof follows by the multiplication with $(1+|\kappa|)^{-b}(1+|\kappa|)^{b}$ directly from the definition of the weighted $L^1$-norm. It is
\begin{align*}
&\|g(\eps\cdot +k_0)D(\cdot)\|_{L_a^1(\R^d)}=\int_{\R^d}(1+|\kappa|)^a|g(\eps\kappa+k_0)||D(\kappa)|\mathrm{d}\kappa = \int_{\R^d}\frac{|g(\eps\kappa+k_0)|}{(1+|\kappa|)^{b}}(1+|\kappa|)^{a+b}|D(\kappa)|\mathrm{d}\kappa\\
&\leq  C\eps^b\sup_{\kappa \in\R^d}\frac{|\kappa|^b}{(1+|\kappa|)^b}\int_{\R^d}(1+|\kappa|)^{a+b}|D(\kappa)|\mathrm{d}\kappa \leq C \eps^b\|D\|_{L^1_{a+b}(\R^d)}.
\end{align*}
\epf
Because
\beqq|\phi(k)|= \bigg|\omega_{n_0}(k)-\sum_{|\alpha|=0}^2\frac{D^\alpha\omega_{n_0}(k_0)}{\alpha !}(k-k_0)^\alpha\bigg|\leq C|k-k_0|^3,\eeqq
we can use Lemma \ref{L:gD-est} with $b=3$ and $a=0$ to obtain
\begin{align}
&\|\phi(\cdot)U_{n_0}^\text{ext}(\cdot,t)\|_{L^1(\B)}= \bigg\|\phi(\cdot)\eps^{1-d}\widehat{A}_r\left(\frac{\cdot-k_0}{\eps},T\right)e^{-\ri\alpha(k,k_0)t}\bigg\|_{L^1(\B)}\notag\\
\leq&  \eps\bigg\|\phi(\eps\cdot+k_0)\widehat{A}_r(\cdot,T)\bigg\|_{L^1\left(\frac{\B-k_0}{\eps}\right)}\leq  C \eps^4 \bigg\|\widehat{A}_r\left(\cdot ,T\right)\bigg\|_{L^1_3\left(\frac{\B-k_0}{\eps}\right)}\leq  C \eps^4 \|\widehat{A}(\cdot,T)\|_{L_3^1(\R^d)}
\label{phiL1Absch}
\end{align}
provided $\widehat{A}(\cdot,T) \in L_3^1(\R^d)$.

At this point we can set our technical parameter $r\in (0,1)$. Estimate \eqref{phiL1Absch} provides the restriction $\beta \geq 3$ and due to \eqref{E:beta-r-rel} we need $r\in (0,2/3)$. We can set, e.g., $r:=\frac{1}{2}.$

Let us now turn to the terms
\beqq \bigg\|\eps^{5-d}\partial_T \widehat{B}_n\left(\frac{\cdot-k_0}{\eps},T\right)\bigg\|_{L^1(\B)}\quad\text{and}\quad \|R_{\eps^{5-d},n}(\cdot,t)\|_{L^1(\B)}\eeqq
for $n\in \N$.  We make use of the asymptotic distribution of the eigenvalues $\omega_n(k)$, see \cite[p.55]{Hoermander3},
\beq\label{asympt}
C_1 n^{2/d}\leq \omega_n(k) \leq C_2 n^{2/d}, \quad k \in \B, n \in \N
\eeq
with $C_1,C_2>0$, in order to pull a $k$-independent factor $c_n(\eps)$ with $(c_n)_n\in l^2_{s/d}$ out of these expressions. Firstly, we have
\beq\label{E:eta-est}
\frac{1}{|\eta_n(k,k_0)|}\leq Cn^{-\frac{2}{d}}, \ n \in \N\setminus\{n_0\}, k \in B_{\eps^r}(k_0)
\eeq
for all $\eps>0$ small enough (see also the discussion below \eqref{Defetan}). Secondly,
\begin{align*}
& |b_{\lambda,\mu,\rho}^n(-k+h,h-l,l,k)|= |\omega_n(k)^{-q}||\langle \sigma(\cdot)\overline{p_{\lambda}}(\cdot,k-h)p_{\mu}(\cdot,h-l)p_{\rho}(\cdot,l),\omega_n(k)^q p_n(\cdot,k)\rangle_{L^2(\P)}|\\
=& |\omega_n(k)^{-q}||\langle \sigma(\cdot)\overline{p_{\lambda}}(\cdot,k-h)p_{\mu}(\cdot,h-l)p_{\rho}(\cdot,l),\mathcal{L}(k)^q p_n(\cdot,k)\rangle_{L^2(\P)}|\\
\leq & |\omega_n(k)^{-q}|\left\|\mathcal{L}(k)^q\left(\sigma(\cdot)\overline{p_{\lambda}}(\cdot,k-h)p_{\mu}(\cdot,h-l)p_{\rho}(\cdot,l)\right)\right\|_{L^2(\P)},
\end{align*}
in which we made use of the self-adjointness of $\mathcal{L}$, the Cauchy-Schwarz-inequality, and the normalization of the Bloch-functions. The asymptotic distribution of eigenvalues in (\ref{asympt}) yields
\beq
|b_{\lambda,\mu,\rho}^n(-k+h,h-l,l,k)| \leq C_{\lambda,\mu,\rho} n^{-\frac{2q}{d}}, \quad n\in \N, k,h,l \in \R^d\label{E:b-est}
\eeq
if $\mathcal{L}(k)^q\left(\sigma(\cdot)\overline{p_{\lambda}}(\cdot,k-h)p_{\mu}(\cdot,h-l)p_{\rho}(\cdot,l)\right)\in L^2(\P)$. Since $\cL$ is of second order, we need to require the $H^{2q}(\P)$-regularity of $\sigma(\cdot)p_{\lambda}(\cdot,k-h)p_{\mu}(\cdot,h-l)p_{\rho}(\cdot,l)$. This holds by the algebra property of $H^{2q}(\P)$ with $q>d/4$ if $\sigma \in H^{2q}(\P)$ and $p_n(\cdot,k)\in H^{2q}(\P)$ for all $n \in \N$. A sufficient condition for the regularity of the Bloch functions $p_n$ is $V \in H_\text{per}^{a}(\P)$, $a>2q+d-2$, cf. Lemma \ref{L:reg}. Estimate \eqref{E:b-est} thus holds if
$$V \in H_\text{per}^{2q+d-2+\delta}(\P), \sigma \in H^{2q}(\P), \ q>d/4, \delta >0.$$

For $\|\eps^{5-d}\partial_T \widehat{B}_n\left(\frac{\cdot-k_0}{\eps},T\right)\|_{L^1(\B)}$ we use \eqref{E:eta-est} and \eqref{E:b-est} with $\lambda=\mu=\rho=n_0$. With the help of Young's inequality for convolutions we get
\begin{align}
&\bigg\|\eps^{5-d}\partial_T \widehat{B}_n\left(\frac{\cdot-k_0}{\eps},T\right)\bigg\|_{L^1(\B)}
= \bigg\|\frac{\eps^{5-d}}{\eta_n(\cdot,k_0)}\int_{B_{2\eps^{r-1}}(0)}\int_{B_{\eps^{r-1}}(0)}\partial_T \left[\widehat{\widebar{A}}_r\left(\frac{\cdot-k_0}{\eps}-\widetilde{h},T\right)\right. \\
& \left.\times  \widehat{A}_r(\widetilde{h}-\widetilde{l},T)\widehat{A}_r(\widetilde{l},T)\right] \widetilde{b}_{n_0,n_0,n_0}^n\left(\tfrac{\cdot-k_0}{\eps},\widetilde{h},\widetilde{l}\right) \mathrm{d}\widetilde{l}\mathrm{d}\widetilde{h}\bigg\|_{L^1(\B)}\notag\\
\leq & C \eps^{5-d}n^{-\frac{2+2q}{d}}\bigg\|\bigg(\partial_T\widehat{\widebar{A}}_r*_{B_{2\eps^{r-1}}(0)}\widehat{A}_r*_{B_{\eps^{r-1}}(0)}\widehat{A}_r+\widehat{\widebar{A}}_r*_{B_{2\eps^{r-1}}(0)}\partial_T\widehat{A}_r*_{B_{\eps^{r-1}}(0)}\widehat{A}_r\notag\\
& \qquad\qquad\qquad+\widehat{\widebar{A}}_r*_{B_{2\eps^{r-1}}(0)}\widehat{A}_r*_{B_{\eps^{r-1}}(0)}\partial_T\widehat{A}_r\bigg)\left(\frac{\cdot-k_0}{\eps},T\right)\bigg\|_{L^1(\B)}\notag\\
\leq & C \eps^{5}n^{-\frac{2+2q}{d}}\|\partial_T\widehat{A}(\cdot,T)\|_{L^1(\R^d)}\|\widehat{A}(\cdot,T)\|^2_{L^1(\R^d)}
\label{AbschAbl}
\end{align}
if $\partial_T\widehat{A}(\cdot,T)\in L^1(\R^d)$, $\widehat{A}(\cdot,T)\in L^1(\R^d),$ $V \in H_\text{per}^{2q+d-2+\delta}(\P), \sigma \in H^{2q}(\P), \ q>d/4, \delta >0$. For the $l^2_{s/d}$-summability in $n$, we have to require $\frac{2s}{d}-\frac{4+4q}{d}<-1$, hence
$$q>\frac{s}{2}+\frac{d}{4}-1.$$

It remains to estimate $\|\vec{R}_{\eps^{5-d}}(\cdot,t)\|_{\cX(s)}$. For this we first introduce a notation for the leading and the higher order parts of the extended ansatz. Namely,
$$
\begin{aligned}
\util^{(0)}_\text{ext}(x,k,t):= U_{n_0}^\text{ext}(k,t)p_{n_0}(x,k), \  \util^{(1)}_\text{ext}(x,k,t) := \sum_{n\neq n_0} U_{n}^\text{ext}(k,t)p_{n}(x,k).
\end{aligned}
$$
We estimate these separately in
\blem\label{L:u0_u1_est}
Let $s>d/2$. If $V\in H_\text{per}^{a}(\P), \sigma \in H_\text{per}^{2q}(\P)$ with $a>\max\{2q+d-2,s+d-2\},$ $q>\tfrac{s}{2}+\tfrac{d}{4}-1$, then
$$\|\util^{(0)}_\text{ext}(\cdot,\cdot,t)\|_{L^1(\B,\Hsper(\P))}\leq C\eps \|\Ahat(\cdot,T)\|_{L^1(\R^d)}, \ \|\util^{(1)}_\text{ext}(\cdot,\cdot,t)\|_{L^1(\B,\Hsper(\P))}\leq C\eps^3 \|\Ahat(\cdot,T)\|^3_{L^1(\R^d)}.$$
\elem
\bpf
For the leading order part we have
$$\|\util^{(0)}_\text{ext}(\cdot,\cdot,t)\|_{L^1(\B,\Hsper(\P))}\leq C\|U^\text{ext}_{n_0}\|_{L^1(\B+k_0)}=C\eps \|\Ahat_r(\cdot,T)\|_{L^1(\eps^{-1}\B)}\leq C\eps \|\Ahat(\cdot,T)\|_{L^1(\R^d)},$$
where the $H^s(\P)$-regularity of $p_{n_0}(\cdot,k_0)$ follows from $V\in H_\text{per}^{a}(\P)$ with $a>s+d-2$, cf. Lemma \ref{L:reg}. 

For $n\neq n_0$ we get
\beq\label{E:Unext-est}
\|U^\text{ext}_{n}(\cdot,t)\|_{L^1(\B+k_0)}=C\eps^3 n^{-\frac{2+2q}{d}}\|\Ahat(\cdot,T)\|^3_{L^1(\R^d)}
\eeq
if $V\in H_\text{per}^{a}(\P), \sigma \in H^{2q}(\P)$ with $a>2q+d-2$ analogously to the calculation leading to \eqref{AbschAbl}. If $q>\tfrac{s}{2}+\tfrac{d}{4}-1$, the $l^2_{s/d}$-summability in $n$ holds and we have
$$\|\util^{(1)}_\text{ext}\|_{L^1(\B,\Hsper(\P))}\leq C\eps^3 \|\Ahat(\cdot,T)\|^3_{L^1(\R^d)}.$$
\epf
We note now that
$$
\begin{aligned}
\vec{R}_{\eps^{5-d}} = & \cD(\sigma \ubartilext^{(1)} *_{\B+k_0} \utilext^{(0)}*_{\B+k_0} \utilext^{(0)}) + 2\cD(\sigma \ubartilext^{(0)} *_{\B+k_0} \utilext^{(1)}*_{\B+k_0} \utilext^{(0)}) \\
&+ 2\cD(\sigma \ubartilext^{(1)} *_{\B+k_0} \utilext^{(1)}*_{\B+k_0} \utilext^{(0)}) +  \cD(\sigma \ubartilext^{(0)} *_{\B+k_0} \utilext^{(1)}*_{\B+k_0} \utilext^{(1)})\\
& +  \cD(\sigma \ubartilext^{(1)} *_{\B+k_0} \utilext^{(1)}*_{\B+k_0} \utilext^{(1)}).
\end{aligned}
$$
With the help of the isomorphism property of $\cD$ and the algebra property in Lemma \ref{L:algeb_L1Hs} we get
$$
\begin{aligned}
\|\vec{R}_{\eps^{5-d}}(\cdot,t)\|_{\cX(s)} &\leq c\|\sigma\|_{H^s(\P)}\left(\|\utilext^{(1)}\|_{L^1(\B,\Hsper(\P))}\|\utilext^{(0)}\|_{L^1(\B,\Hsper(\P))}^2+\|\utilext^{(1)}\|_{L^1(\B,\Hsper(\P))}^2\|\utilext^{(0)}\|_{L^1(\B,\Hsper(\P))}\right.\\
& \qquad \left.+\|\utilext^{(1)}\|_{L^1(\B,\Hsper(\P))}^3\right)\\
&\leq c \eps^5,
\end{aligned}
$$
where $c<\infty$ if $\Ahat(\cdot,T)\in L^1(\R^d)$. In the last step we used Lemma \ref{L:u0_u1_est}.

To summarize, we have for any $s>d/2$
\beq \|\vec{Res}(\cdot,t)\|_{\cX(s)}\leq C_\text{res} \eps^4\label{ResAbsch}\eeq
if
$$\widehat{A}(\cdot,T)\in L_3^1(\R^d), \ \partial_T\widehat{A}(\cdot,T)\in L^1(\R^d), \ \sigma\in H_\text{per}^{2q}(\P), \text{ and } V \in H_\text{per}^{a}(\P)$$
for some
$$a>\max\{2q+d-2,s+d-2\}, \ q>\max\left\{\frac{s}{2}+\frac{d}{4}-1,\frac{d}{4}\right\}.$$

\subsection{Estimation of the error}\label{S:esterror}

We first rewrite equation \eqref{E:ODE-Un} in the vector form
\beq \ri\partial_t\vec{U}=W(k)\vec{U}+\vec{F}(\vec{U},\vec{U},\vec{U})\label{GleichungfuerU}\eeq
with
\beqq W(k):=\text{diag}((\omega_n(k))_{n\in\N}), \quad \vec{U}(k,t):=(U_n(k,t))_{n\in\N},\eeqq
\beqq\vec{F}(\vec{U},\vec{U},\vec{U}):=\left(\langle \sigma(\cdot)(\widetilde{\widebar{u}}*_{\B}\widetilde{u}*_{\B}\widetilde{u})(\cdot,k,t),p_n(\cdot,k)\rangle_{L^2(\P)}\right)_{n\in\N}.\eeqq
Analogously to $\Uvec$ we define $\Uvecext$ 
and consider first the error generated by the extended ansatz, i.e.
$$\vec{E}(k,t):=\vec{U}(k,t)-\vec{U}^\text{ext}(k,t).$$
Equation \eqref{GleichungfuerU} is equivalent to
\begin{align*}
\ri\partial_t\vec{E}&=W(k)\vec{E}+G(\vec{U}^\text{ext},\vec{E}),
\end{align*}
where
\beqq
G(\vec{U}^\text{ext},\vec{E}):=(Res_n(k,t))_{n\in\N}+\vec{F}(\vec{U},\vec{U},\vec{U})-\vec{F}(\vec{U}^\text{ext},\vec{U}^\text{ext},\vec{U}^\text{ext}).
\eeqq
The variation of constants produces
\beq \label{E:int-eq-E}
\vec{E}(t)=\vec{E}(0)+\int_0^t S(t-\tau)G(\vec{U}^\text{ext},\vec{E})(\tau)\mathrm{d}\tau
\eeq
in which $S(t)_{t\geq 0}$, given by $S(t)=e^{-\ri tW}:\cX(s)\rightarrow\cX(s)$, is a uniformly bounded strongly continuous semigroup. We first need an estimate of the nonlinearity $\|G(\vec{U}^\text{ext},\vec{E})(\cdot,t)\|_{\cX(s)}$. For the sake of brevity we leave out the $t-$dependence in the following calculation. Due to \eqref{ResAbsch} it is
\beqq \|G(\vec{U}^\text{ext},\vec{E})\|_{\cX(s)}\leq C_\text{res} \eps^4 + \|\vec{F}(\vec{U},\vec{U},\vec{U})-\vec{F}(\vec{U}^\text{ext},\vec{U}^\text{ext},\vec{U}^\text{ext})\|_{\cX(s)}.\eeqq
Next, using the isomorphism property of $\cD$ and the algebra property in Lemma \ref{L:algeb_L1Hs}, we have
\begin{align*}
&\left\|\vec{F}(\vec{U},\vec{U},\vec{U})-\vec{F}(\vec{U}^\text{ext},\vec{U}^\text{ext},\vec{U}^\text{ext})\right\|_{\cX(s)} \leq C \left\|\sigma\left(\widetilde{\widebar{u}}*_{\B}\widetilde{u}*_{\B}\widetilde{u}-\widetilde{\widebar{u}}_\text{ext}*_{\B}\widetilde{u}_\text{ext}*_{\B}\widetilde{u}_\text{ext}\right)\right\|_{L^1(\B,\Hsper(\P))}\\
&\hspace{2cm}\leq C\left(\|\widetilde{u}_\text{ext}\|^2_{L^1(\B,\Hsper(\P))}\|\widetilde{e}\|_{L^1(\B,\Hsper(\P))} + \|\widetilde{u}_\text{ext}\|_{L^1(\B,\Hsper(\P))}\|\widetilde{e}\|^2_{L^1(\B,\Hsper(\P))}+\|\widetilde{e}\|^3_{L^1(\B,\Hsper(\P))}\right)\\
&\hspace{2cm}\leq C\left(\|\Uvecext\|^2_{\cX(s)}\|\vec{E}\|_{\cX(s)} + \|\Uvecext\|_{\cX(s)}\|\vec{E}\|^2_{\cX(s)} + \|\vec{E}\|_{\cX(s)}^3 \right),
\end{align*}
where $\widetilde{e}(x,k,t):=\sum_{n\in \N}E_n(k,t)p_n(x,k)$.

Using Lemma \ref{L:u0_u1_est}, we have
$$\|\Uvecext(\cdot,t)\|_{\cX(s)}\leq c (\|\utilext^{(0)}(\cdot,\cdot,t)\|_{L^1(\B,\Hsper(\P))}+\|\utilext^{(1)}(\cdot,\cdot,t)\|_{L^1(\B,\Hsper(\P))}) \leq c\eps$$
if $\Ahat(\cdot,T)\in L^1(\R^d)$ with $a>\max\{2q+d-2,s+d-2\},$ $q>\tfrac{s}{2}+\tfrac{d}{4}-1$. In that case we thus have
$$
\left\|G(\vec{U}^\text{ext},\vec{E})\right\|_{\cX(s)} \leq C_\text{res} \eps^4 + c_1 \left(\eps^2\|\vec{E}\|_{\cX(s)}+\eps\|\vec{E}\|^2_{\cX(s)} +\|\vec{E}\|^3_{\cX(s)}\right).
$$
Next, we need to estimate $\|\vec{E}(\cdot,0)\|_{\cX(s)}$ in \eqref{E:int-eq-E}. Because of our assumption $\Uvec(\cdot,0)=\Uvecapp(\cdot,0)$, it suffices to estimate $\Uvecext-\Uvecapp$. This is done in Lemma \ref{L:Uapp-Uext} and leads to
\beq
\|\vec{E}(\cdot,0)\|_{\cX(s)} \leq c_0 \eps^2 \label{AbschE(0)}
\eeq
if $V\in H_\text{per}^{a}(\P), \sigma \in H_\text{per}^{2q}(\P)$ with some $a>2q+d-2$ and $q>\frac{s}{2}+\frac{d}{4}$ and if $\Ahat(\cdot,0)\in L^1_\beta(\R^d)$ with some $\beta>2q+d$.
The above estimates and the unitary property of $S$ imply
\beqq
\|\vec{E}(\cdot,t)\|_{\cX(s)}\leq c_0\eps^2 +tC_{Res} \eps^4  + c_1 \int_0^t \eps^2\|\vec{E}(\cdot,\tau)\|_{\cX(s)}+\eps \|\vec{E}(\cdot,\tau)\|^2_{\cX(s)} +\|\vec{E}(\cdot,\tau)\|^3_{\cX(s)} \mathrm{d}\tau.
\eeqq
For convenience, we define $\vec{\cE}:=\eps^{-2}\vec{E}$ and using a bootstrapping argument and Gronwall's inequality, we show that if $\Ahat(\cdot,T)\in L^1_\beta(\R^d)$ for all $T\in [0,T_0]$, then there is $M>0$ and $\eps_0>0$ such that $\|\vec{\cE}(\cdot,t)\|_{\cX(s)}\leq M$ for all $t\in [0,\eps^{-2}T_0]$ if $\eps\in (0,\eps_0)$.

As long as $\|\vec{\cE}(\cdot,t)\|_{\cX(s)}\leq M$ (which, by continuity, holds on a non-empty time interval for any $M>c_0$), we have
$$\|\vec{\cE}(\cdot,t)\|_{\cX(s)}\leq c_0 + t\eps^2\left(C_\text{res} + c_1(\eps M^2 + \eps^2M^3)\right) + c_1\eps^2\int_0^t\|\vec{\cE}(\cdot,\tau)\|_{\cX(s)}\dd \tau.$$
Gronwall's inequality then produces
$$
\|\vec{\cE}(\cdot,t)\|_{\cX(s)}\leq \left[c_0 + t\eps^2\left(C_\text{res} + c_1(\eps M^2 + \eps^2M^3)\right)\right]e^{c_1\eps^2 t}
$$
as long as $\|\vec{\cE}(\cdot,t)\|_{\cX(s)}\leq M$. Next, we make a suitable choice of $M$ and $\eps_0$.
Namely, we set
$$M:=\left[c_0 + T_0 (C_\text{res} +1)\right]e^{c_1T_0}$$
and choose $\eps_0$ so small that
$$c_1(\eps_0 M^2 + \eps_0^2M^3)\leq 1.$$
Then $\|\vec{\cE}(\cdot,t)\|_{\cX(s)}\leq M$ for all $t\in [0,\eps^{-2}T_0]$ if $\eps\in (0,\eps_0)$, as desired.

Finally, the estimate in Theorem \ref{T:main} follows by applying Lemma \ref{L:sup_control}, the isomorphism property of $\cD$, Lemma \ref{L:Uapp-Uext}, as well as the triangle inequality.
\begin{align*}
\sup_{x\in \R^d}\sup_{t\in[0,T_0\eps^{-2}]}|u(x,t)-u_\text{app}(x,t)|
\leq & \sup_{x\in \R^d}\sup_{t\in[0,T_0\eps^{-2}]}C \|\widetilde{u}(x,\cdot,t)-\widetilde{u}_\text{app}(x,\cdot,t)\|_{L^1(\B,H_\text{per}^{s}(\P))}\\
\leq &C\sup_{t\in[0,T_0\eps^{-2}]} \left(\|(\vec{U}-\vec{U}^\text{ext})(\cdot,t)\|_{\cX(s)}+ \|(\vec{U}^\text{ext}-\vec{U}^\text{app})(\cdot,t)\|_{\cX(s)}\right)\\
\leq &C\sup_{t\in[0,T_0\eps^{-2}]} \left(\|\vec{E}(\cdot,t)\|_{\cX(s)}+\|(\vec{U}^\text{ext}-\vec{U}^\text{app})(\cdot,t)\|_{\cX(s)}\right)\leq C \eps^2
\end{align*}
if $\eps\in (0,\eps_0)$. In addition, due to $\widetilde{u}-\widetilde{u}_\text{app} \in L^1(\B,\Hsper(\P))$ Lemma \ref{L:sup_control} gives also the decay $|u(x,t)-u_\text{app}(x,t)|\rightarrow 0$ for $|x|\rightarrow\infty$. 

\blem\label{L:Uapp-Uext}
If $V\in H_\text{per}^{a}(\P)$ and $\sigma \in H_\text{per}^{2q}(\P)$ with some $a>2q+d-2$, $q>\frac{s}{2}+\frac{d}{4}$, and $\Ahat(\cdot,T)\in L^1_\beta(\R^d)$ with some $\beta>2q+d$, then there is $c>0$ such that
$$\|(\Uvecapp-\Uvecext)(\cdot,t)\|_{\cX(s)}\leq c\left(\eps^2 \|\Ahat(\cdot,T)\|_{L^1_\beta(\R^d)} + \eps^3\|\Ahat(\cdot,T)\|^3_{L^1(\R^d)} \right)$$
for all $\eps>0$ small enough.
\elem
\bpf
We estimate first $\|(\Uapp_{n_0}-\Uext_{n_0})(\cdot,t)\|_{L^1(\B+k_0)}$ and then separately $\|\Uext_n(\cdot,t)\|_{L^1(\B+k_0)}$ and $\|\Uapp_n(\cdot,t)\|_{L^1(\B+k_0)}$ with $n\neq n_0$. Recall that $\Uapp_n$ and $\Uext_n$ are defined in \eqref{KoeffappAns} and \eqref{KoefferwAnsneu} respectively.
With the notation
$$\pi_{n,m}(k):=\langle p_{n_0}(\cdot,k_0-m),p_n(\cdot,k)\rangle_{L^2(\P)}$$
we have
\begin{align*}
\|(\Uext_{n_0}&-\Uapp_{n_0})(\cdot,t)\|_{L^1(\B+k_0)} \\
&=\eps^{1-d}\left\|\cX_{B_{\eps^{1/2}}(k_0)}(\cdot)\Ahat\left(\frac{\cdot-k_0}{\eps},T\right)-\sum_{m\in \Z^d}\Ahat\left(\frac{\cdot-k_0+m}{\eps},T\right)\pi_{n_0,m}(\cdot)e^{-\ri m\cdot v_g t}\right\|_{L^1(\B+k_0)}\\
&\leq \eps^{1-d}\left(\left\|(1-\pi_{n_0,0}(\cdot))\Ahat\left(\frac{\cdot-k_0}{\eps},T\right)\right\|_{L^1(B_{\eps^{1/2}}(k_0))}+\left\|\Ahat\left(\frac{\cdot-k_0}{\eps},T\right)\right\|_{L^1(\R^d\setminus B_{\eps^{1/2}}(k_0))}\right).
\end{align*}
Recall that we have set $r:=1/2$ above.
Using the normalization of the Bloch waves and their Lipschitz continuity with respect to $k$, we get
$$|1-\pi_{n_0,0}(k)|=|\langle p_{n_0}(\cdot,k_0),p_{n_0}(\cdot,k_0)-p_{n_0}(\cdot,k)\rangle_{L^2(\P)}|\leq L |k-k_0|$$
for all $k\in \B+k_0$. Next, the tail of $\Ahat$ is estimated via
$$\left\|\Ahat\left(\frac{\cdot-k_0}{\eps},T\right)\right\|_{L^1(\R^d\setminus B_{\eps^{1/2}}(k_0))} =\eps^d \|\widetilde{A}(\cdot,T)\|_{L^1(\R^d)}\leq \eps^{d+\beta/2}\|\Ahat(\cdot,T)\|_{L^1_\beta(\R^d)}$$
for any $\beta>0$ as shown in \eqref{L1beta}. Hence
\beq\label{E:n0-diff-est}
\|(\Uext_{n_0}-\Uapp_{n_0})(\cdot,t)\|_{L^1(\B+k_0)} \leq \eps^2L\|\Ahat(\cdot,T)\|_{L^1_1(\R^d)} + c\eps^{\beta/2+1}\|\Ahat(\cdot,T)\|_{L^1_\beta(\R^d)}.
\eeq

Next, we consider the indices $n\neq n_0$. $\|\Uext_n(\cdot,t)\|_{L^1(\B+k_0)}$ is estimated in \eqref{E:Unext-est}. For
$$\|\Uapp_n(\cdot,t)\|_{L^1(\B+k_0)}=\left\|\eps^{1-d}\sum_{m\in \Z^d}\Ahat\left(\frac{\cdot-k_0+m}{\eps},T\right)\pi_{n,m}(\cdot)\right\|_{L^1(\B+k_0)}$$
we first estimate $|\pi_{n,m}(k)|$. For $n\neq n_0$, $q\in \N$, and $V\in H^{a}(\P), a>2q+d-3$, is
\begin{align}
|\pi_{n,0}(k)|&=|\omega_n(k)|^{-q}|\langle p_{n_0}(\cdot,k_0)-p_{n_0}(\cdot,k),\omega_n(k)^qp_n(\cdot,k)\rangle_{L^2(\P)}|\notag\\
&\leq cn^{-\frac{2q}{d}}\|\cL(k)^q(p_{n_0}(\cdot,k_0)-p_{n_0}(\cdot,k))\|_{L^2(\P)}\notag\\
&\leq c(\|p_{n_0}(\cdot,k_0)\|_{H^{2q-1}(\P)})n^{-\frac{2q}{d}} |k-k_0|, \quad k\in \B+k_0 \label{E:pi-n-0-est}
\end{align}
thanks to the Lipschitz continuity of the Bloch functions in Lemma \ref{L:Lip}.

For $n\neq n_0$ and $m\neq 0$ we estimate 
\begin{align*}
|\pi_{n,m}(k)|&\leq cn^{-\frac{2q}{d}}\|\cL(k)^q(p_{n_0}(\cdot,k_0-m)-p_{n_0}(\cdot,k))\|_{L^2(\P)} \\
& \leq cn^{-\frac{2q}{d}}\left(\omega_{n_0}(k)^q+\|\cL(k)^q p_{n_0}(\cdot,k-m)\|_{L^2(\P)}\right).
\end{align*}

Because $p_{n_0}(x,k-m)=p_{n_0}(x,k)e^{\ri m\cdot x}$, it is $\|\cL(k)^q p_{n_0}(\cdot,k-m)\|_{L^2(\P)}\leq c |m|^{2q}$ if $V\in H^{a}(\P)$, $a>2q+d-2$ (see Lemma \ref{L:reg}) and we get
\beq\label{E:pi-n-m-est}
|\pi_{n,m}(k)|\leq cn^{-\frac{2q}{d}}|m|^{2q}\quad k\in \B+k_0, m \in \Z^d, n \neq n_0.
\eeq
We write next
\begin{align*}
\|\Uapp_n(\cdot,t)\|_{L^1(\B+k_0)}\leq &\eps^{1-d}\left\|\Ahat\left(\frac{\cdot}{\eps},T\right) \pi_{n,0}(\cdot+k_0)\right\|_{L^1(\B)} \\
&+ \eps^{1-d}\sum_{m\in \Z^d\setminus\{0\}}\left\|\Ahat\left(\frac{\cdot}{\eps},T\right)\pi_{n,m}(\cdot+k_0-m)\right\|_{L^1(\B+m)}.
\end{align*}
In the second term we have got only the tail of $\Ahat$, which we estimate via
\begin{align}\label{E:Atail-m-est}
\int_{\B+m}\left|\Ahat\left(\frac{k}{\eps},T\right)\right|\dd k&=\eps^d \int_{\eps^{-1}(\B+m)}\left|\Ahat(\kappa,T)\right|\dd \kappa \notag\\
&\leq \eps^d\|\Ahat(\cdot,T)\|_{L^1_\beta(\R^d)}\sup_{\kappa\in \eps^{-1}(\B+m)}(1+|\kappa|)^{-\beta} \leq c\eps^{d+\beta}|m|^{-\beta}\|\Ahat(\cdot,T)\|_{L^1_\beta(\R^d)}.
\end{align}
Combining \eqref{E:pi-n-0-est}, \eqref{E:pi-n-m-est}, and \eqref{E:Atail-m-est} yields
\beq\label{E:Uappn-est}
\|\Uapp_n(\cdot,t)\|_{L^1(\B+k_0)}\leq cn^{-\frac{2q}{d}}\left(\eps^2 \|\Ahat(\cdot,T)\|_{L^1_1(\R^d)} + \eps^{1+\beta} \|\Ahat(\cdot,T)\|_{L^1_\beta(\R^d)} \sum_{m\in \Z^d\setminus\{0\}}|m|^{2q-\beta}\right), \ n\neq n_0.
\eeq
For the $m-$sum we need to require $\beta>2q+d$, i.e. we need $\Ahat(\cdot,T)\in L^1_\beta(\R^d)$ with $\beta>2q+d$.

Finally, if $q>\frac{s}{2}+\frac{d}{4}$, then \eqref{E:n0-diff-est}, \eqref{E:Unext-est} and \eqref{E:Uappn-est} prove the result.
\epf

The required regularity on $V$ and $\sigma$ follows from the conditions in \eqref{ResAbsch}, in Lemma \ref{L:Uapp-Uext}, and from $s>d/2$ needed in Lemma \ref{L:sup_control}. Namely, it suffices to set $q>\tfrac{d}{2}$  because then there exists an $s>d/2$ such that $q>\tfrac{s}{2}+\tfrac{d}{4}$. The requirements in \eqref{ResAbsch} and Lemma \ref{L:Uapp-Uext} produce the conditions
$$\sigma\in H_\text{per}^{d+\delta}(\P), V \in H_\text{per}^{2d-2+\delta}(\P) \text{ with some } \delta>0.$$

\section{A Nonlinear Wave Equation}\label{S:NLW}

We consider next the semilinear wave equation (NLW)
\beq\label{E:NLW}
\pa_t^2 u=\chi_1(x)\Delta u-\chi_2(x)u-\chi_3(x)u^3, \ (x,t)\in \R^d\times \R
\eeq
with $\chi_j:\R\to \R$, $\chi_j(x+2\pi e_i)=\chi_j(x)$ for all $x\in\R^d$, $1\leq j\leq
3$ and $1\leq i\leq d$
and $\chi_1(x)\geq \gamma >0, \ \chi_2(x) >0$  for all $x\in \R^d.$

An obvious change compared to the GP is that now the equation is real, hence the
solution ansatz is to be chosen real. In practice, for our
approximate wavepacket ansatz, we take twice the real part of the
GP-ansatz (with $p_{n_0}$ being eigenfunctions of a different eigenvalue problem). The analysis is for the most part analogous to the GP case. The main modification is caused by the fact that the nonlinearity $u^3$
(as opposed to $|u|^2u$) generates higher harmonics. For an ansatz at a fixed frequency $\omega_0$ the nonlinearity $u^3$ generates residual terms also at frequencies $\pm 3 \omega_0$, which have to be accounted for.

In order to keep the notation simple, we recycle and change the
meaning of some of the symbols from Section \ref{S:GP}. We also work on
the same periodicity cell $\P=(0,2\pi]^d$.

\subsection{Introduction, first order formulation and the approximate ansatz}

In the wave equation case temporal frequencies are square roots of the positive eigenvalues $\lambda_n(k)$ of $\mathcal{L}(k)$. Hence we define
$$\omega_n:=\sqrt{\lambda_n}, \quad \omega_{-n}:=-\sqrt{\lambda_n}, \quad n \in \N.$$
The eigenfunctions $p_n(x,k)$ are chosen such that for each $k\in \B$ they form an $L^2_{\chi_1}(\P)$-orthonormal system.

Applying $\cD\cT$, i.e. the Bloch transformation and the
$L^2_{\chi_1}(\P)$-projection onto the Bloch eigenfunctions $p_n(\cdot,k)$, we get
\beq\label{E:ODE_Un2}
\partial_t^2 U_n(k,t)=-\omega_n^2(k)U_n(k,t)-\langle \chi_3(\cdot)(\widetilde{u}*_\B \widetilde{u}*_\B \widetilde{u})(\cdot,k,t),p_n(\cdot,k)\rangle_{L^2_{\chi_1}}, \ k\in \B, t>0,
\eeq 
where $\widetilde{u}(x,k,t)=\sum_{n\in \N} U_n(k,t)p_n(x,k)$. We reformulate equation \eqref{E:ODE_Un2} as a first order ODE-system using similar variables to those in \cite{BSTU06}. Namely, we define
$V_n(k,t):=\omega_n(k)^{-1}\partial_t U_n(k,t)$ and $Z_n:=
(U_n,V_n)^T$. This yields
$$ \partial_t Z_n(k,t)=\bpm 0 & \omega_n(k)\\-\omega_n(k)& 0\epm Z_n(k,t)- \bpm 0 \\ \omega_n(k)^{-1}\langle \chi_3(\cdot)(\widetilde{u}*_\B \widetilde{u}*_\B\widetilde{u})(\cdot,k,t),p_n(\cdot,k)\rangle_{L^2_{\chi_1}}\epm, $$
which after diagonalization is equivalent to
$$ \ri\partial_t \widetilde{Z}_n(k,t)=\bpm \omega_n(k) & 0 \\ 0 & \omega_{-n}(k) \epm \widetilde{Z}_n(k,t)-\ri \mathbb{Q}\bpm 0 \\ \omega_n(k)^{-1}\langle \chi_3(\cdot)(\widetilde{u}*_\B \widetilde{u}*_\B\widetilde{u})(\cdot,k,t),p_n(\cdot,k)\rangle_{L^2_{\chi_1}}\epm$$
for
$$\widetilde{Z}_n:=\mathbb{Q} Z_n, \quad \mathbb{Q}:=
\frac{1}{\sqrt{2}}\bspm 1 & \ri \\ 1 & -\ri\espm.$$
Finally, we write
$$\bspm \widetilde{z}_n\\ \widetilde{z}_{-n}\espm:=
\widetilde{Z}_n.$$ 
The vector $(\tilde{z}_n)_{n\in \Z_0}$, where $\Z_0:=\Z\setminus\{0\}$, satisfies the system
\beq\label{E:ODE_zn} 
\ri\partial_t \widetilde{z}_n(k,t)-\omega_n(k)\widetilde{z}_n(k,t)-\frac{1}{\sqrt{2}}\omega_n(k)^{-1}\langle
\chi_3(\cdot)(\widetilde{u}*_\B
\widetilde{u}*_\B\widetilde{u})(\cdot,k,t),p_n(\cdot,k)\rangle_{L^2_{\chi_1}}=0
\eeq
for all $n\in \Z_0$.

Let us define $\widetilde{v}(x,k,t):=\sum_{n\in \N} V_n(k,t)p_n(x,k)$. For each $k\in \B$ we denote by $\widetilde{\cD}(k)$ the
expansion-diagonalization operator, i.e.
$$\widetilde{\cD}(k):\bspm
\widetilde{u}(\cdot,k,t)\\\widetilde{v}(\cdot,k,t)\espm \mapsto
\left(\widetilde{Z}_n(k,t)\right)_{n\in\N}.$$

Formally, the inverse is given by
 \beq\label{E:IsomInv}
\widetilde{\cD}(k)^{-1}\left(\left(\widetilde{z}_n(k,t)\right)_{n\in\Z\setminus\{0\}}\right)(x)=\frac{1}{\sqrt{2}}\bpm
\sum_{n\in\N}(\widetilde{z}_n+\widetilde{z}_{-n})(k,t)p_n(x,k) \\
\sum_{n\in\N}\ri(\widetilde{z}_{-n}-\widetilde{z}_n)(k,t)p_n(x,k)\epm=\bpm
\widetilde{u}(x,k,t)\\\widetilde{v}(x,k,t)\epm, \eeq where
$\widetilde{v}(x,k,t)=\sum_{n\in\N}V_n(k,t)p_n(x,k)$. We denote then
by $\widetilde{\cD}$ the collective operator
$\widetilde{\cD}:\bspm \widetilde{u}(\cdot,\cdot,t)\\ \widetilde{v}(\cdot,\cdot,t)\espm
\mapsto \left(\widetilde{Z}_n(\cdot,t)\right)_{n\in\N}.$

Similarly to Lemma \ref{L:Diag} we have
\blem\label{L:isom-NLW} 
For $s>\frac{d}{2}$ the map
$$ \widetilde{\cD}:L^1(\B,\Hsper(\P)\times \Hsper(\P))\rightarrow \cX(s):=L^1(\B,l^2_{s/d}(\N)\times l^2_{s/d}(\N))\cong L^1(\B,l^2_{s/d}(\Z\setminus\{0\}))$$
with
$$ \|\left(\widetilde{z}_n(\cdot,t)\right)_{n\in\Z\setminus\{0\}}\|_{\cX(s)}=\int_\B \|\left(\widetilde{z}_n(k,t)\right)_{n\in\Z\setminus\{0\}}\|_{l^2_{s/d}}\mathrm{d}k=\int_\B \left(\sum_{n\in\Z\setminus\{0\}}|n|^{2s/d}|\widetilde{z}_n(k,t)|^2\right)^{1/2}\mathrm{d}k$$
is an isomorphism.
\elem

We define next our approximate ansatz for a pulse solution. For a
given $k_0 \in \B $ and $n_0\in \N$ we define again $v_g:=\nabla
\omega_{n_0}(k_0)$ and look for a solution $u$ of \eqref{E:NLW}	
close to the ansatz 
\beq\label{E:uapp2} u_{\text{app}}(x,t)=\eps
A(\eps(x-v_g t),\eps^2 t)p_{n_0}(x,k_0)e^{\ri(k_0\cdot x -\omega_0t)}+
\text{c.c.}, \eeq 
where  c.c. denotes the complex conjugate. Note that this ansatz is analogous to \eqref{E:uapp} except we take here the real part. In the Bloch variables $\util=\cT u$ the ansatz reads
\beq\label{E:uapp2_Bloch} 
\begin{aligned}
\tilde{u}_{\text{app}}(x,k,t)=&\eps^{1-d}
\sum_{m\in \Z^d}\Ahat\left(\frac{k-k_0+m}{\eps},T\right)e^{-\ri\Lambda_1(k+m,k_0)t}p_{n_0}(x,k_0-m)\\
+ &\eps^{1-d}
\sum_{m\in \Z^d}\widehat{\widebar{A}}\left(\frac{k+k_0+m}{\eps},T\right)e^{-\ri\Lambda_{-1}(k+m,k_0)t}p_{n_0}(x,-k_0-m),
\end{aligned}
\eeq
where
$$\Lambda_j(k,k_0):=j\omega_0 + (k-jk_0)\cdot v_g, \quad T:=\eps^2 t.$$

Our main result for the nonlinear wave equation is
\bthm\label{T:main-NLW}
Let $\chi_{1,2}\in H_\text{per}^{2d-2+\delta}(\P)$ and $\chi_{3} \in H_\text{per}^{d+\delta}(\P)$ with some $\delta>0$ satisfy \eqref{E:per-gen} and \eqref{E:ellip}. Assume that the Bloch eigenvalue $\lambda_{n_0}(k)$ of $\cL(k)$ in \eqref{E:Bloch-gen} is simple at $k=k_0$ with the corresponding eigenfunction $p_{n_0}(\cdot,k_0)$. 
In addition assume that the nonresonance condition 
\beq\label{E:NR-cond}
\inf_{\stackrel{(n,j)\in \{\Z_0\times \{\pm 1,\pm 3\}\}}{(n,j)\notin\{(n_0,1),(-n_0,-1)\}}}\left|j\omega_{n_0}(k_0)-\omega_{n}(jk_0)\right|>0
\eeq
holds.
Consider the ansatz $u_{\text{app}}$ in \eqref{E:uapp2}. For every solution $A$ of \eqref{E:NLS} with $\nu$ given in \eqref{E:nu} and
with the regularity $\widehat{A}\in C([0,T_0]; L^2(\R^d)\cap L^1_{\beta+1}(\R^d))$, $\pa_T\Ahat \in C([0,T_0];L^1_{\beta}(\R^d))$ for some $T_0>0$ and $\beta>2d$ there exist $\eps_0>0$ and $C>0$ such that for all $\eps\in (0,\eps_0)$ the solution of \eqref{E:NLW} with initial data $u(x,0)=u_\text{app}(x,0)$ is continuous and satisfies
$$\sup_{x\in \R^d}|u(x,t)-u_\text{app}(x,t)|\leq C \eps^2 \quad \text{for all }t\in [0,T_0\eps^{-2}]$$
as well as the decay $u(x,t)\to 0$ for $|x|\to \infty$.
\ethm

Because we work again with the expansion coefficients $(U_n)_{n\in \N}$ of $\tilde{u}$, we calculate these next for $\tilde{u}_{app}$. They are
$$
\begin{aligned}
\Uapp_n(k,t)=&\eps^{1-d}\sum_{m\in \Z^d}\Ahat\left(\frac{k-k_0+m}{\eps},T\right) e^{-\ri \Lambda_1(k+m,k_0)t}\pi^+_{n,m}(k)\\
+ &\eps^{1-d}\sum_{m\in \Z^d}\Ahatbar\left(\frac{k+k_0+m}{\eps},T\right) e^{-\ri \Lambda_{-1}(k+m,k_0)t}\pi^-_{n,m}(k)
\end{aligned}
$$
with the notation
$$\pi^\pm_{n,m}(k):=\langle p_{n_0}(\cdot,\pm k_0-m),p_n(\cdot,k)\rangle_{L^2_{\chi_1}(\P)}.$$

For $\tilde{z}_n^{\text{app}}(k,t):=\frac{1}{\sqrt{2}}\left(\Uapp_n(k,t)+\ri \omega_n^{-1}(k)\pa_t\Uapp_n(k,t)\right), n \in \N$ we get
$$\tilde{z}_n^{\text{app}} = \tilde{z}_n^{\text{app,lead}} + R^{\text{app}}_n,$$
where
$$
\begin{aligned}
\tilde{z}_n^{\text{app,lead}}(k,t):=&\eps^{1-d}\sqrt{2}\sum_{m\in \Z^d}\Ahat\left(\frac{k-k_0+m}{\eps},T\right) e^{-\ri \Lambda_1(k+m,k_0)t}\pi^+_{n,m}(k),\\
R^{\text{app}}_n(k,t):=&\eps^{1-d}\sum_{m\in \Z^d}\frac{\omega_{n_0}(k_0)+(k+m-k_0)\cdot v_g-\omega_n(k)}{\sqrt{2}\omega_n(k)}\times \\
 &\quad \times \Ahat\left(\frac{k-k_0+m}{\eps},T\right)e^{-\ri \Lambda_1(k+m,k_0)t}\pi^+_{n,m}(k)\\
 - &\eps^{1-d}\sum_{m\in \Z^d}\frac{\omega_{n_0}(k_0)-(k+m+k_0)\cdot v_g-\omega_n(k)}{\sqrt{2}\omega_n(k)}\times \\
&\quad \times \Ahatbar\left(\frac{k+k_0+m}{\eps},T\right)e^{-\ri \Lambda_{-1}(k+m,k_0)t}\pi^-_{n,m}(k)\\
 + &\eps^{3-d}\frac{\ri}{\sqrt{2}\omega_n(k)}\sum_{m\in \Z^d}\pa_T\Ahat\left(\frac{k-k_0+m}{\eps},T\right)e^{-\ri \Lambda_{1}(k+m,k_0)t}\pi^+_{n,m}(k)\\
+  &\eps^{3-d}\frac{\ri}{\sqrt{2}\omega_n(k)}\sum_{m\in \Z^d}\pa_T\Ahatbar\left(\frac{k+k_0+m}{\eps},T\right)e^{-\ri \Lambda_{-1}(k+m,k_0)t}\pi^-_{n,m}(k).
\end{aligned}
$$
And for $\tilde{z}_{-n}^{\text{app}}(k,t):=\frac{1}{\sqrt{2}}\left(\Uapp_n(k,t)-\ri \omega_n^{-1}(k)\pa_t\Uapp_n(k,t)\right), n \in \N$ we get
$$
\tilde{z}_{-n}^{\text{app}}(k,t) = \sqrt{2}\eps^{1-d}\sum_{m\in \Z^d}\Ahatbar\left(\frac{k+k_0+m}{\eps},T\right) e^{-\ri \Lambda_{-1}(k+m,k_0)t}\pi^-_{n,m}(k) - R^{\text{app}}_{n}(k,t).
$$

\subsection{Definition of the extended ansatz}
Analogously to the GP case the approximate ansatz does not produce a small enough residual. Hence, we define a modified (``extended'') ansatz. An added complexity compared to the GP case is that the residual evaluated at the approximate ansatz includes higher harmonics generated by the nonlinearity. We get $u_\text{app}^3(x,t) = \eps^3(A^3p_{n_0}(x,k_0)^3e^{3\ri (k_0\cdot x-\omega_0 t)}+3|A|^2A|p_{n_0}(x,k_0)|^2p_{n_0}(x,k_0)e^{\ri (k_0\cdot x-\omega_0 t)}+\text{c.c.})$. These higher harmonics are accounted for by correction terms in the extended ansatz.
\begin{align}\label{E:defzn0ext}
\widetilde{z}_{n_0}^\text{ext}(k,t):=\eps^{1-d}\sqrt{2}&\sum_{m\in\Z^d}\Ahat_r\left(\frac{k+m-k_0}{\eps},T\right)e^{-\ri \Lambda_1(k+m,k_0)t}\\
&+\eps^{3-d}\sum_{j\in\{-1,\pm 3\}}\sum_{m\in\Z^d}\Ahat_{n_0,j}\left(\frac{k+m-jk_0}{\eps},T\right)e^{-\ri \Lambda_j(k+m,k_0)t}, \  (k,t)\in \B\times\R\notag
\end{align}
\begin{align}\label{E:defz-n0ext}
\widetilde{z}_{-n_0}^\text{ext}(k,t):=\eps^{1-d}\sqrt{2}&\sum_{m\in\Z^d}\Ahatbar_r\left(\frac{k+m+k_0}{\eps},T\right)e^{-\ri \Lambda_{-1}(k+m,k_0)t}\\
&+\eps^{3-d}\sum_{j\in\{1,\pm 3\}}\sum_{m\in\Z^d}\Ahat_{-n_0,j}\left(\frac{k+m-jk_0}{\eps},T\right)e^{-\ri \Lambda_j(k+m,k_0)t}, \  (k,t)\in \B\times\R\notag
\end{align}
Further for $n\neq \pm n_0$
\beq\label{E:defznext}
\widetilde{z}_n^\text{ext}(k,t):=\eps^{3-d}\sum_{j\in\{\pm 1, \pm 3\}}\sum_{m\in\Z^d}\Ahat_{n,j}\left(\frac{k+m-jk_0}{\eps},T\right)e^{-\ri \Lambda_j(k+m,k_0)t}, (k,t)\in \B\times\R.
\eeq
All envelope functions $\Ahat_r, \Ahat_{n,j}$ are chosen to have  a compact support, namely
$$\text{supp}(\Ahat_r(\cdot,T))\subset B_{\eps^{r-1}}(0), \ \text{supp}(\Ahat_{n,j}(\cdot,T))\subset B_{3\eps^{r-1}}(0) \ \text{for all } n,j \ \text{and some } r\in (0,1).$$
Once again, $\Ahat_r$ will be chosen as the cut-off of the Fourier  transform of a solution to the effective NLS equation and $\Ahat_{n,j}$ will be explicit functions of $\Ahat_r$.

Note that in contrary to the GP case we cannot drop the $m$-sums in the formulas for $\widetilde{z}_{n}^\text{ext}(k,t), k \in \B$ despite the compact supports of the envelopes because $\widetilde{z}_{n}^\text{ext}(k,t)$ consists of several terms centered at different $k$-points. Outside $\B$ the functions $\widetilde{z}_n^\text{ext}(\cdot,t)$ are defined by $1$-periodicity in each coordinate direction.

For the variable $\widetilde{u}$ the extended ansatz is given by
\beq\label{E:extans}
\widetilde{u}_\text{ext}(x,k,t):=\frac{1}{\sqrt{2}}\sum_{n\in\N}(\widetilde{z}_n^\text{ext}+\widetilde{z}_{-n}^\text{ext})(k,t)p_n(x,k).
\eeq

The residual corresponding to equation \eqref{E:ODE_zn} is
\beq\label{E:res2}
Res_n(k,t):=\ri \pa_t \widetilde{z}_n^\text{ext}(k,t)-\omega_n(k)\widetilde{z}_n^\text{ext}(k,t)-s_{n}(k,t), \quad n \in \Z_0,
\eeq
where 
\beq\label{E:s_n}
s_{n}(k,t):=\frac{1}{\sqrt{2}\omega_n(k)}\langle \chi_3(\cdot)(\widetilde{u}_\text{ext}*_\B \widetilde{u}_\text{ext}*_\B\widetilde{u}_\text{ext})(\cdot,k,t),p_n(\cdot,k)\rangle_{L^2_{\chi_1}}.
\eeq
The functions $\Ahat_r$ and $\Ahat_{n,j}$ will be chosen in Sec. \ref{S:res-calc-NLW} such that the residual is sufficiently small and such that the leading order part
$$\widetilde{z}_{n_0}^\text{ext, lead}(k,t):=\eps^{1-d}\sqrt{2}\sum_{m\in\Z^d}\Ahat_r\left(\frac{k+m-k_0}{\eps},T\right)e^{-\ri \Lambda_1(k+m,k_0)t}$$
is close to $\widetilde{z}_{n_0}^\text{app}(k,t)$. This and some direct estimates will lead to the smallness of $\widetilde{z}^\text{app}-\widetilde{z}^\text{ext}$ in Lemma \ref{L:Zapp-Zext} similarly to 
Lemma \ref{L:Uapp-Uext} for the GP. Note that in the leading order part of the extended ansatz the $m-$sum can be dropped due to the support of $\Ahat_r$. We have namely
$$\widetilde{z}_{n_0}^\text{ext, lead}(k,t)=\eps^{1-d}\sqrt{2}\Ahat_r\left(\frac{k-k_0}{\eps},T\right)e^{-\ri \Lambda_1(k,k_0)t}, \ k \in \B+k_0.$$
Analogously, we define $\widetilde{z}_{-n_0}^\text{ext, lead}$ as the $O(\eps^{1-d})$-term in  $\widetilde{z}_{-n_0}^\text{ext}$ and we set $\widetilde{z}_{n}^\text{ext, lead}:=0$ for $|n|\neq n_0.$

We also use the notation
$$
\begin{aligned}
\widetilde{u}_\text{ext}^{(0)}(x,k,t)&:=\frac{1}{\sqrt{2}}\sum_{n\in\N}(\widetilde{z}_n^\text{ext,lead}+\widetilde{z}_{-n}^\text{ext,lead})(k,t)p_n(x,k),\\
\widetilde{u}_\text{ext}^{(1)}&:=\widetilde{u}_\text{ext}-\widetilde{u}_\text{ext}^{(0)}.
\end{aligned}
$$

\subsection{Calculation of the residual}\label{S:res-calc-NLW}

The difficult part of the residual is the nonlinearity $s_{n}$ in \eqref{E:s_n}. We split the nonlinearity into the leading order part and the higher order rest
$$s_{n}=s^\text{lead}_{n}+s^\text{hot}_{n},$$
where
$$s^\text{lead}_{n}(k,t):=\frac{1}{\sqrt{2}\omega_n(k)}\langle \chi_3(\cdot)(\widetilde{u}_\text{ext}^{(0)}*_\B \widetilde{u}_\text{ext}^{(0)}*_\B\widetilde{u}_\text{ext}^{(0)})(\cdot,k,t),p_{|n|}(\cdot,k)\rangle_{L^2_{\chi_1}}$$
and
$$
\begin{aligned}
s^\text{hot}_{n}(k,t):=\frac{1}{\sqrt{2}\omega_n(k)}\Big[&3\langle \chi_3(\cdot)(\widetilde{u}_\text{ext}^{(0)}*_\B \widetilde{u}_\text{ext}^{(0)}*_\B\widetilde{u}_\text{ext}^{(1)})(\cdot,k,t),p_{|n|}(\cdot,k)\rangle_{L^2_{\chi_1}} \\ 
+&3\langle \chi_3(\cdot)(\widetilde{u}_\text{ext}^{(0)}*_\B \widetilde{u}_\text{ext}^{(1)}*_\B\widetilde{u}_\text{ext}^{(1)})(\cdot,k,t),p_{|n|}(\cdot,k)\rangle_{L^2_{\chi_1}}\\
+ & \left.\langle \chi_3(\cdot)(\widetilde{u}_\text{ext}^{(1)}*_\B \widetilde{u}_\text{ext}^{(1)}*_\B\widetilde{u}_\text{ext}^{(1)})(\cdot,k,t),p_{|n|}(\cdot,k)\rangle_{L^2_{\chi_1}}\right].
\end{aligned}
$$
A similar calculation to that in Section \ref{S:res-calc} leads to
$$s^\text{lead}_{n}=\sum_{J\in \{\pm 1,\pm 3\}}  s^\text{lead}_{n,J},
$$
where for each $J$ the term $s^\text{lead}_{n,J}$ denotes those contributions which are concentrated near $Jk_0$, i.e.,
$$
\begin{aligned}
s^\text{lead}_{n,J}(k,t):=\eps^{3-d}&\sum_{\stackrel{j_1,j_2,j_3\in \{\pm n_0\}}{S_{j_1,j_2,j_3}=J}}\int_{B_{2\eps^{r-1}}(0)}\int_{B_{\eps^{r-1}}(0)} \Gamma_{j_1}(\Ahat_r)(\kappa-h,T)\Gamma_{j_2}(\Ahat_r)(h-l,T)\Gamma_{j_3}(\Ahat_r)(l,T)\times\\
&\times \beta^n_{j_1,j_2,j_3}(\kappa-h,h-l,l,\kappa) \dd l\dd h e^{-\ri \Lambda_{J}(Jk_0+\eps \kappa,k_0)t}, \ k \in \B  + Jk_0,
\end{aligned}
$$
where 
$$
\begin{aligned}
&S_{j_1,j_2,j_3}:=\sum_{m=1}^3\sign(j_m), \kappa:=\frac{k-S_{j_1,j_2,j_3}k_0}{\eps}, \Gamma_j(\Ahat_r) := \begin{cases} \Ahat_r, \text{ if } j=n_0,\\
\Ahatbar_r, \text{ if } j=-n_0, \end{cases}\\
& \text{ and } \\
&\beta^n_{j_1,j_2,j_3}(\kappa-h,h-l,l,\kappa):= \frac{2}{\omega_n(k)}\langle \chi_3(\cdot)p_{|j_1|}(\cdot,-\sign(j_1)k_0+\eps(h-\kappa))\times \\
&\hspace{2cm} \times p_{|j_2|}(\cdot,\sign(j_2)k_0+\eps(h-l))p_{|j_3|}(\cdot,\sign(j_3)k_0+\eps l),p_{|n|}(\cdot,S_{j_1,j_2,j_3}k_0+\eps \kappa)\rangle_{L^2_{\chi_1}}.
\end{aligned}
$$

The residual has the leading order part at the formal order $O(\eps^{3-d})$. We call this part $Res^\text{lead}$. Below, in Sec. \ref{S:resleadn0k0}-\ref{S:resleadremain} we first study $Res^\text{lead}$ and then in Sec. \ref{S:reshot} the higher order rest $Res^\text{hot}:=Res-Res^\text{lead}$.
We first consider $Res_{n_0}(k,t)$ for $k$ near $k_0$ in Sec. \ref{S:resleadn0k0} and $Res_{-n_0}(k,t)$ near $k=-k_0$ in Sec. \ref{S:resleadmn0mk0}. The envelope $\Ahat_r$ will have to be chosen as the (truncated) Fourier transform of a solution to an effective NLS equation in order to make these residual components small. Afterwards, in Sec. \ref{S:reshot}, the other contributions to $Res(k,t)$ will be considered and $\Ahat_{n,j}$ chosen in dependence on $\Ahat_r$ in order to make $Res_n$ small.

\subsubsection{The leading order part $Res^\text{lead}_{n_0}(k,t)$ for $k\in B_{3\eps^{r}}(k_0)$}\label{S:resleadn0k0}

The relevant leading order nonlinear term of $\chi_{B_{3\eps^{r}}(k_0)}(k)Res_{n_0}(k,t)$ is that part of $s_{n_0}^\text{lead}$ which is concentrated near $k_0$. It is  the term
$$
\begin{aligned}
s^\text{lead}_{n_0,1}(k,t)=&3\eps^{3-d}\int_{B_{2\eps^{r-1}}(0)}\int_{B_{\eps^{r-1}}(0)} \Ahatbar_r(\kappa-h,T)\Ahat_r(h-l,T)\Ahat_r(l,T)\times\\
&\quad \times\beta^{n_0}_{-n_0,n_0,n_0}(h-\kappa,h-l,l,\kappa) \dd l\dd h e^{-\ri \Lambda_{1}(k_0+\eps \kappa,k_0)t},\\
\end{aligned}
$$
where
\beq\label{E:kappa}
\kappa := \frac{k-k_0}{\eps}.
\eeq
Just like in the GP case we approximate the integrand by a double convolution of $\Ahat$ on $\R^d$ and the coefficient $\beta^{n_0}_{-n_0,n_0,n_0}(\kappa-h,h-l,l,\kappa)$ by
\beq\label{E:nu}
\nu:=\frac{3}{\sqrt{2}}\beta^{n_0}_{-n_0,n_0,n_0}(0,0,0,0).
\eeq
This leads to
$$
s^\text{lead}_{n_0,1}(k,t)=
 \eps^{3-d}\sqrt{2}\nu \chi_{B_{\eps^{r}}(0)}(k)(\Ahatbar*_{\R^d}\Ahat*_{\R^d}\Ahat)(\kappa,T)e^{-\ri \Lambda_{1}(k_0+\eps \kappa,k_0)t} + \psi(k,t),
$$
where
$$
\begin{aligned}
\psi(k,t):=& \eps^{3-d}\nu\left[(\Ahatbar_r *_{B_{2\eps^{r-1}}(0)}\Ahat_r *_{B_{\eps^{r-1}}(0)} \Ahat_r)(\kappa,T) - \chi_{B_{\eps^{r-1}}(0)}(\kappa)(\Ahatbar *_{\R^d} \Ahat *_{\R^d} \Ahat)(\kappa,T) \right]\times \\
&\quad \times e^{-\ri \Lambda_{1}(k_0+\eps \kappa,k_0)t}\\
&+ \eps^{3-d}\int_{B_{2\eps^{r-1}}(0)}\int_{B_{\eps^{r-1}}(0)} \Ahatbar_r(\kappa-h,T)\Ahat_r(h-l,T)\Ahat_r(l,T) \times\\
&\quad \times (3\beta^{n_0}_{-n_0,n_0,n_0}(h-\kappa,h-l,l,\kappa)-\sqrt{2}\nu) \dd l\dd h e^{-\ri \Lambda_{1}(k_0+\eps \kappa,k_0)t}.
\end{aligned}
$$

Next, the linear terms in $Res_{n_0}(k,t)$ near $k_0$ are precisely $(\ri \pa_t-\omega_{n_0}(k))\tilde{z}^\text{ext,lead}(k,t)$, i.e.
$$
\begin{aligned}
&\sqrt{2}\left[\eps^{1-d}(\Lambda_1(k_0+\eps\kappa,k_0)-\omega_{n_0}(k))\Ahat_r(\kappa,T) + \ri \eps^{3-d}\pa_T\Ahat_r(\kappa,T)\right]e^{-\ri \Lambda_1(k_0+\eps \kappa,k_0)t}\\
&=\sqrt{2}\left[\eps^{3-d}\left(\ri \pa_T\Ahat_r(\kappa,T)-\tfrac{1}{2}\kappa^TD^2\omega_{n_0}(k_0)\kappa\Ahat_r(\kappa,T)\right)+\eps^{1-d}\phi(k)\Ahat_r(\kappa,T)\right]e^{-\ri \Lambda_1(k_0+\eps \kappa,k_0)t},
\end{aligned}
$$
where $|\phi(k)|\leq c |k-k_0|^3$  for some $c>0$ and all $k\in \B$ as follows by the Taylor expansion of $\omega_{n_0}(k)$ in $k=k_0$.

Taking the sum of this linear term with $s^\text{lead}_{n_0,1}$, we see that the leading order part of the residual $\chi_{B_{3\eps^{r}}(k_0)}(k)Res_{n_0}(k,t)$ reduces to
$$
\chi_{B_{3\eps^{r}}(k_0)}(k)Res_{n_0}^\text{lead}(k,t)=\eps^{1-d}\sqrt{2}\phi(k)\Ahat_r(\kappa,T)e^{-\ri \Lambda_1(k_0+\eps \kappa,k_0)t} + \psi(k,t)
$$
provided we choose again
\beq\label{E:Ahatr}
\Ahat_r(\kappa,T):=\chi_{B_{\eps^{r-1}}}(\kappa)\Ahat(\kappa,T),
\eeq
where $A$ is a solution of \eqref{E:NLS} with $\nu$ given in \eqref{E:nu}.

\subsubsection{The leading order part $Res^\text{lead}_{-n_0}(k,t)$ for $k\in B_{3\eps^{r}}(-k_0)$}\label{S:resleadmn0mk0}

The same calculation as above produces
$$
\begin{aligned}
\chi_{B_{3\eps^{r}}(-k_0)}(k)Res^\text{lead}_{-n_0}(k,t)=\chi_{B_{3\eps^{r}}(k_0)}(-k)\overline{Res^\text{lead}_{n_0}}(-k,t)
\end{aligned}
$$

\subsubsection{The remaining components of the leading order residual}
\label{S:resleadremain}

The rest of the residual at $O(\eps^{3-d})$ includes firstly the nonlinearity terms $s_{n,j}^\text{lead}$ not accounted for in Sec. \ref{S:resleadn0k0}-\ref{S:resleadremain}, i.e. all $(n,j) \in \Z_0\times \{\pm 1,\pm 3\}$ other than $(n_0,1)$ and $(-n_0,-1)$. Secondly, it includes the linear terms evaluated at the correction $\widetilde{z}_{n}^\text{ext}-\widetilde{z}_{n}^\text{ext, lead}$. Hence, we have 
$$
\begin{aligned}
\chi_{B_{3\eps^r}(jk_0)}(k) Res^\text{lead}_{n}(k,t) = \eps^{3-d}&\left[(j\omega_0 +(k-jk_0)^T v_g-\omega_n(k)) \Ahat_{n,j}\left(\frac{k-jk_0}{\eps},T\right)e^{-\ri \Lambda_j(k,k_0)t}\right.\\
&\quad \left. - s^\text{lead}_{n,j}(k,t)\right],
\end{aligned}
$$
for 
$$(n,j) \in \{\Z_0\times \{\pm 1,\pm 3\}\} \setminus \{(n_0,1), (-n_0,-1)\}.$$
This part of the residual is eliminated by the choice
\beq\label{E:Ahatnj}
\Ahat_{n,j}\left(\frac{k-jk_0}{\eps},T\right):=\frac{s^\text{lead}_{n,j}(k,t)e^{\ri \Lambda_j(k,k_0)t}}{j\omega_0 +(k-jk_0)^T v_g-\omega_n(k)}.
\eeq
The above denominators stay bounded away from zero for all $k \in B_{\eps^r}(jk_0)$ with $\eps>0$ small enough as long as the non-resonance condition \eqref{E:NR-cond} holds.

\subsubsection{The higher oder part of the residual}
\label{S:reshot}

The rest of the residual is the higher order nonlinearity $\eps^{5-d}s_{n}^\text{ext,hot}(k,t), n \in \Z_0$, and the time derivative applied to the $\eps^2 t$-dependence of $\Ahat_{n,j}$. In detail
$$
\begin{aligned}
Res_{n_0}^\text{hot}(k,t)=&\eps^{5-d}\left(\sum_{m\in \Z^d}\sum_{j\in \{-1,\pm 3\}}\ri\pa_T\Ahat_{n_0,j}\left(\frac{k+m-jk_0}{\eps},T\right)e^{-\ri \Lambda_j(k+m,k_0)t}+s_{n_0}^\text{ext,hot}(k,t)\right),\\
Res_{-n_0}^\text{hot}(k,t)=&\eps^{5-d}\left(\sum_{m\in \Z^d}\sum_{j\in \{1,\pm 3\}}\ri\pa_T\Ahat_{-n_0,j}\left(\frac{k+m-jk_0}{\eps},T\right)e^{-\ri \Lambda_j(k+m,k_0)t}+s_{-n_0}^\text{ext,hot}(k,t)\right),\\
Res_{n}^\text{hot}(k,t)=&\eps^{5-d}\left(\sum_{m\in \Z^d}\sum_{j\in \{\pm 1,\pm 3\}}\ri\pa_T\Ahat_{n,j}\left(\frac{k+m-jk_0}{\eps},T\right)e^{-\ri \Lambda_j(k+m,k_0)t}+s_{n}^\text{ext,hot}(k,t)\right),\\
\end{aligned}
$$
$n \in \Z_0\setminus\{\pm n_0\}$.

\subsection{Estimation of the residual}

After the choice of $\Ahat_r$ in \eqref{E:Ahatr} and of $\Ahat_{n,j}$ in  \eqref{E:Ahatnj} the residual reduces to 
$$
\begin{aligned}
Res_{n_0}(k,t)&= \eps^{1-d}\sqrt{2}\phi(k)\Ahat_r(\tfrac{k-k_0}{\eps},T)e^{-\ri \Lambda_1(k,k_0)t} + \psi(k,t) + Res_{n_0}^\text{hot}(k,t),\\
Res_{-n_0}(k,t)&=\overline{Res_{n_0}}(-k,t),\\
Res_{n}(k,t)&=Res_{n}^\text{hot}(k,t), \ n\in \Z_0\setminus\{\pm n_0\}.
\end{aligned}
$$
The estimates of individual terms are analogous to the GP case in Sec. \ref{S:res-est}. Therefore, we mostly present only the results.

Just like in \eqref{E:psi-est} and \eqref{phiL1Absch} we get
\beq\label{E:phiest-NLW}
\|\psi(\cdot,t)\|_{L^1(\B)}\leq C \eps^4 \|\Ahat(\cdot,T)\|^3_{L^1_\beta(\R^d)} \text{  and  } \eps^{1-d}\|\phi(\cdot)\Ahat_r\left(\tfrac{\cdot-k_0}{\eps},T\right)\|_{L^1(\B)}\leq C\eps^4 \|\Ahat(\cdot,T)\|_{L^1_3(\R^d)}
\eeq
if $\beta\geq \tfrac{1}{1-r}$ and $\Ahat(\cdot,T)\in L^1_3(\R^d)\cap L^1_\beta(\R^d)$. Because $r:=1/2$ can be selected, it suffices to choose  $\Ahat(\cdot,T)\in L^1_3(\R^d)$.

Also for the higher order terms the estimates are analogous. Note that $\omega_n$ is now the square root of the Bloch eigenvalue, such that the asymptotic distribution is 
\beq\label{asympt-NLW}
C_1 n^{1/d}\leq \omega_n(k) \leq C_2 n^{1/d}, \quad k \in \B, n \in \N.
\eeq
Nevertheless, we obtain the same decay of $Res_n$ in $n$ as in the GP case because of the $\frac{1}{\omega_n(k)}$ factor in $s^\text{lead}_{n}$ and $s^\text{hot}_{n}$. In detail,
$$|\beta_{j_1,j_2,j_3}^n(q,r,s,t)|\leq cn^{-\frac{2q+1}{d}}$$
if $\chi_3,p_{j_1},p_{j_2},p_{j_3}\in H_\text{per}^{2q}(\P)$ with $q>d/4$. By Lemma \ref{L:reg} this regularity of the Bloch functions follows if 
$$\chi_1,\chi_2\in H_\text{per}^{a}(\P), \quad a>2q+d-2, q>\frac{d}{4}.$$

This leads to 
$$\eps^{5-d}\left\|\sum_{m\in \Z^d}\sum_{j\in \{\pm 1,\pm 3\}}\pa_T\Ahat_{n,j}\left(\frac{\cdot+m-jk_0}{\eps},T\right)\right\|_{L^1(\B)}\leq c\eps^5n^{-\frac{2q+2}{d}}\|\pa_T\Ahat(\cdot,T)\|_{L^1(\R^d)}\|\Ahat(\cdot,T)\|^2_{L^1(\R^d)}.$$

Next, for $s_{n}^\text{ext,hot}$ we first get the same estimates on $\widetilde{u}_\text{ext}^{(0)}$ and $\widetilde{u}_\text{ext}^{(1)}$ as in Lemma 
\ref{L:u0_u1_est} if $\chi_1,\chi_2\in H^{a}_\text{per}(\P), \chi_3\in H^{2q}(\P)$ with $a>\max\{2q+d-2,s+d-2\},$ $q>\tfrac{s}{2}+\tfrac{d}{4}-1$.

For $\vec{\tilde{s}}^\text{ext,hot}$ with $\tilde{s}_{n}^\text{ext,hot}(k,t):=\frac{1}{\omega_n(k)}s_{n}^\text{ext,hot}(k,t)$, we get
$\|\vec{\tilde{s}}^\text{ext,hot}(\cdot,t)\|_{\cX(s)}\leq c\eps^5$
with $c=c(\|\Ahat(\cdot,T)\|_{L^1(\R^d)})$. Hence also 
$$\|\vec{\tilde{s}}^\text{ext,hot}(\cdot,t)\|_{\cX(s)}\leq c\eps^5$$
if $\Ahat(\cdot,T)\in L^1(\R^d)$.

In summary, we have again for any $s>d/2$
\beq 
\|\vec{Res}(\cdot,t)\|_{\cX(s)}\leq C_\text{res} \eps^4\label{ResAbsch-NLW}\eeq
if
$$\widehat{A}(\cdot,T)\in L_3^1(\R^d), \ \partial_T\widehat{A}(\cdot,T)\in L^1(\R^d), \ \chi_1,\chi_2\in H_\text{per}^{a}(\P), \text{ and } \chi_3\in H_\text{per}^{2q}(\P)$$
for some
$$a> \max\{2q+d-2,s+d-2\}, \ q>\max\left\{\frac{s}{2}+\frac{d}{4}-1,\frac{d}{4}\right\}.$$

\subsection{Estimation of the error}

Restricting to real solutions $u$ is equivalent to assuming $\tilde{z}_n(k,t)=\overline{\tilde{z}}_{-n}(-k,t)$ for all $n \in \N, k \in \B$.
This restriction allows us to consider system \eqref{E:ODE_zn} only for $n\in \N$, where $\util(x,k,t)=\tfrac{1}{\sqrt{2}}\sum_{n\in \N}\left(\tilde{z}_n(k,t)+\overline{\tilde{z}}_{n}(-k,t)\right)p_n(x,k)$. Similarly to \eqref{GleichungfuerU} we write this system in the vector form
\beq 
\ri\partial_t\vec{Z}=W(k)\vec{Z}+\vec{F}(\vec{Z},\vec{Z},\vec{Z}),
\label{GleichungfuerZ}
\eeq
where
\beqq 
W(k):=\text{diag}((\omega_n(k))_{n\in\N}), \quad \vec{Z}(k,t):=(\tilde{z}_n(k,t))_{n\in\N},
\eeqq
\beqq\vec{F}(\vec{Z},\vec{Z},\vec{Z}):=\left(\langle \chi_3(\cdot)(\widetilde{u}*_{\B}\widetilde{u}*_{\B}\widetilde{u})(\cdot,k,t),p_n(\cdot,k)\rangle_{L_{\chi_1}^2(\P)}\right)_{n\in\N}.\eeqq
The error generated by the extended ansatz $\vec{E}:=\vec{Z}-\vec{Z}^{\text{ext}}$, where $\vec{Z}^{\text{ext}}=(\tilde{z}_n^\text{ext})_n$, satisfies again 
\begin{align*}
\ri\partial_t\vec{E}&=W(k)\vec{E}+G(\vec{Z}^\text{ext},\vec{E}),
\end{align*}
where
\beqq
G(\vec{Z}^\text{ext},\vec{E}):=(Res_n(k,t))_{n\in\N}+\vec{F}(\vec{Z},\vec{Z},\vec{Z})-\vec{F}(\vec{Z}^\text{ext},\vec{Z}^\text{ext},\vec{Z}^\text{ext}).
\eeqq

The analogous result to Lemma \ref{L:u0_u1_est} holds, i.e. 
if $\chi_1,\chi_2\in H_\text{per}^{a}(\P), \chi_3 \in H_\text{per}^{2q}(\P)$ with $a>\max\{2q+d-2,s+d-2\},$ $q>\tfrac{s}{2}+\tfrac{d}{4}-1$ and $s>d/2$, then
\beq\label{E:util01_est}
\|\util^{(0)}_\text{ext}(\cdot,\cdot,t)\|_{L^1(\B,\Hsper(\P))}\leq c\eps \|\Ahat(\cdot,T)\|_{L^1(\R^d)}, \ \|\util^{(1)}_\text{ext}(\cdot,\cdot,t)\|_{L^1(\B,\Hsper(\P))}\leq c\eps^3 \|\Ahat(\cdot,T)\|^3_{L^1(\R^d)}.
\eeq
Just like in Sec. \ref{S:esterror} Gronwall's inequality then produces 
$$\|\vec{E}(\cdot,t)\|_{\cX(s)}\leq c\eps^2 \quad \text{for all } t\in [0, T_0\eps^{-2}]$$
if $\Ahat(\cdot,T)\in L^1_\beta(\R^d)$ for all $T\in [0,T_0]$, $\beta>2q+d$ and $\|\vec{E}(\cdot,0)\|_{\cX(s)} \leq c \eps^2$. This last assumption holds again thanks to the estimate of $\|(\vec{Z}^\text{ext}-\vec{Z}^\text{app})(\cdot,t)\|_{\cX(s)}$ in Lemma \ref{L:Zapp-Zext}, where 
$$\vec{Z}^\text{app}(k,t):=(\tilde{z}^\text{app}_n(k,t))_{n\in\N},$$ 
and due the fact that $u(x,0)=u_\text{app}(x,0)$.

Hence, it remains to estimate the difference $\vec{Z}^\text{ext}-\vec{Z}^\text{app}$. The process is similar to Lemma \ref{L:Uapp-Uext} but we have to estimate also the terms $R^\text{app}_n$ in $\tilde{z}^\text{app}_n$. In the rest of the proof recall that $r=1/2$.

\blem\label{L:Zapp-Zext}
If $\chi_1,\chi_2\in H_\text{per}^{2q-2}(\P)$ and $\chi_3 \in H_\text{per}^{2q}(\P)$ with some $q>\frac{s}{2}+\frac{d}{4}$, and $\Ahat(\cdot,T)\in L^1_{\beta+1}(\R^d)$ and $\pa_T\Ahat(\cdot,T)\in L^1_\beta(\R^d)$ with some $\beta>2q+d$, then there is a constant $c>0$ such that
$$\|(\vec{Z}^\text{ext}-\vec{Z}^\text{app})(\cdot,t)\|_{\cX(s)}\leq c\left(\eps^2 \|\Ahat(\cdot,T)\|_{L^1_{\beta+1}(\R^d)} + \eps^3\|\Ahat(\cdot,T)\|^3_{L^1(\R^d)}+\eps^4\|\pa_T\Ahat(\cdot,T)\|_{L^1_\beta(\R^d)}\right)$$
for all $\eps>0$ small enough.
\elem
\bpf
We get just like in \eqref{E:n0-diff-est} the estimate
$$\|(\tilde{z}_{n_0}^\text{app}-\tilde{z}^{\text{ext,lead}}_{n_0})(\cdot,t)\|_{L^1(\B)} \leq \eps^2L\|\Ahat(\cdot,T)\|_{L^1_1(\R^d)} + c\eps^{\beta/2+1}\|\Ahat(\cdot,T)\|_{L^1_\beta(\R^d)}
$$
for some $L,c>0$. Next, analogously to \eqref{E:Uappn-est} we get for $n\neq n_0$
$$\|\tilde{z}_n^\text{app,lead}(\cdot,t)\|_{L^1(\B)}\leq cn^{-\frac{2q}{d}}\left(\eps^2\|\Ahat(\cdot,T)\|_{L^1_1(\R^d)}+\eps^{1+\beta}\|\Ahat(\cdot,T)\|_{L^1_\beta(\R^d)}\sum_{m\in \Z^d\setminus\{0\}}|m|^{2q-\beta} \right)$$
if $\beta>2q+d$.
We can thus conclude
\begin{align}
\|(\vec{Z}^{\text{app,lead}}-\vec{Z}^\text{ext})(\cdot,t)\|_{\cX(s)}\leq &\|\util^{(1)}_\text{ext}(\cdot,\cdot,t)\|_{L^1(\B,\Hsper(\P))}+\|(\tilde{z}_n^\text{app,lead})_{n\in \N\setminus\{n_0\}}\|_{\cX(s)} \notag\\
& + \|(\tilde{z}_{n_0}^\text{app,lead}-\tilde{z}^\text{ext,lead}_{n_0})(\cdot,t)\|_{L^1(\B)}\notag\\
\leq & c\left(\eps^2 \|\Ahat(\cdot,T)\|_{L^1_1(\R^d)} + \eps^{\beta/2+1}\|\Ahat(\cdot,T)\|_{L^1_\beta(\R^d)}+\eps^3\|\Ahat(\cdot,T)\|_{L^1(\R^d)}^3\right)\label{E:Zapplead-Zext}
\end{align}
if $q>\tfrac{s}{2}+\tfrac{d}{4}$, where $\util^{(1)}_\text{ext}$ is estimated in \eqref{E:util01_est}.

It remains to estimate $\|(R^\text{app}_n(\cdot,t))_n\|_{\cX(s)}$. In $R^\text{app}_{n_0}$ we note that $\omega_{n_0}(k_0)+(k+m-k_0)^T v_g-\omega_{n_0}(k)$ is quadratic in $k+m-k_0$ such that 
$$
\left\|\sum_{m\in \Z^d}(\omega_{n_0}(k_0)+(\cdot+m-k_0)^T v_g-\omega_{n_0}(\cdot))\Ahat\left(\frac{\cdot-k_0+m}{\eps},T\right)\right\|_{L^1(\B)}\leq \eps^{2+d} \|\Ahat(\cdot,T)\|_{L^1_2(\R^d)}.
$$
As a result
$$\|R^\text{app}_{n_0}(\cdot,t)\|_{L^1(\B)}\leq c\eps^3(\|\Ahat(\cdot,T)\|_{L^1_2(\R^d)}+\|\pa_T\Ahat(\cdot,T)\|_{L^1(\R^d)}).$$
For $n\neq n_0$ we first have
$$
\begin{aligned}
&\left\|\sum_{m\in \Z^d}\frac{\omega_{n_0}(k_0)+(\cdot+m-k_0)^T v_g-\omega_{n}(\cdot)}{\omega_n(\cdot)}\Ahat\left(\frac{\cdot-k_0+m}{\eps},T\right)\pi_{n,m}^+(\cdot)\right\|_{L^1(\B)}\\
&\qquad \leq c\left\|\sum_{m\in \Z^d}(1+|\cdot+m-k_0|)\Ahat\left(\frac{\cdot-k_0+m}{\eps},T\right)\pi_{n,m}^+(\cdot)\right\|_{L^1(\B)}\\
&\qquad \leq c \left\|(1+|\cdot-k_0|)\Ahat\left(\frac{\cdot-k_0}{\eps},T\right)\pi_{n,0}^+(\cdot)\right\|_{L^1(\B)} + c\sum_{m\in \Z^d\setminus\{0\}}\left\|\Ahat\left(\frac{\cdot-k_0}{\eps},T\right)\pi_{n,0}^+(\cdot)\right\|_{L^1(\B+m)} \\
&\qquad + c\sum_{m\in \Z^d\setminus\{0\}}\left\| |\cdot-k_0|\Ahat\left(\frac{\cdot-k_0}{\eps},T\right)\pi_{n,0}^+(\cdot)\right\|_{L^1(\B+m)}\\
&\qquad \leq c \eps^{d}n^{-\frac{2q}{d}}\left[\eps\|\Ahat(\cdot,T)\|_{L^1_2(\R^d)}+ \|\Ahat(\cdot,T)\|_{L^1_\beta(\R^d)}\sum_{m\in \Z^d\setminus\{0\}}|m|^{2q}\sup_{\kappa\in \eps^{-1}(\B+m)}(1+|\kappa|)^{-\beta}\right. \\
&\qquad \left. +\eps\|\Ahat(\cdot,T)\|_{L^1_{\beta+1}(\R^d)}\sum_{m\in \Z^d}(|m|^{2q}+1)\sup_{\kappa\in \eps^{-1}(\B+m)}(1+|\kappa|)^{-\beta}\right]\\
&\qquad \leq c \eps^{d}n^{-\frac{2q}{d}}\left[\eps\|\Ahat(\cdot,T)\|_{L^1_2(\R^d)} + \eps^\beta\|\Ahat(\cdot,T)\|_{L^1_\beta(\R^d)} + \eps^{\beta+1}\|\Ahat(\cdot,T)\|_{L^1_{\beta+1}(\R^d)}\right]
\end{aligned}
$$
for all $\beta > 2q+d.$ Here we have used, firstly, the bounds \eqref{E:pi-n-0-est} and \eqref{E:pi-n-m-est}, which hold also for $\pi_{n,m}^\pm$ and, secondly, the fact that $\pi_{n,m}^\pm(k)=\pi_{n,0}^\pm(k+m)$.

Similarly, we have for $n\neq n_0$ and the $\pa_T\Ahat$ term in $R^\text{app}_n$
$$
\begin{aligned}
\left\|\sum_{m\in \Z^d}\pa_T\Ahat\left(\frac{\cdot-k_0+m}{\eps},T\right)\pi_{n,m}^+(\cdot)\right\|_{L^1(\B)} \leq c n^{-\frac{2q}{d}}\left(\eps^{d+1}\|\pa_T\Ahat(\cdot,T)\|_{L^1_1(\R^d)}+\eps^{d+\beta}\|\pa_T\Ahat(\cdot,T)\|_{L^1_\beta(\R^d)}\right)
\end{aligned}
$$
if $\beta>2q+d$. The terms with $\Ahatbar$ and $\pa_T\Ahatbar$ are treated analogously, where we note that $\omega_{n_0}(k_0)-(k+m+k_0)^T v_g - \omega_n(k)=\omega_{n_0}(-k_0)+(k+m+k_0)^T (\nabla_k\omega_{n_0}(-k_0)) - \omega_n(k)$.
This leads to the estimate
$$\|(R^\text{app}_n(\cdot,t))_n\|_{\cX(s)} \leq c \left(\eps^2 \|\Ahat(\cdot,T)\|_{L^1_{\beta+1}(\R^d)}+\eps^4\|\pa_T\Ahat(\cdot,T)\|_{L^1_{\beta}(\R^d)}\right)$$
if $q>\frac{s}{2}+\frac{d}{4}$. We have used here $\beta>2$, which follows from $\beta>2q+d>s+\tfrac{3}{2}d>2d$. Together with \eqref{E:Zapplead-Zext} this proves the statement.
\epf

By the same argument as at the end of Sec. \ref{S:esterror} we conclude that the regularity $\chi_{1},\chi_{2}\in H_\text{per}^{2d-2+\delta}(\P)$ and $\chi_{3} \in H_\text{per}^{d+\delta}(\P)$ with some $\delta>0$ is sufficient for the requirements in Lemma \ref{L:Zapp-Zext} and in the residual estimate \eqref{ResAbsch-NLW}. Therefore, Theorem \ref{T:main-NLW} is proved.

\section{Numerical Example}\label{S:num}

We present a numerical example of an approximate solitary wave in the two dimensional ($d=2$) GP \eqref{E:GP} as predicted by Theorem \ref{T:main} with a Townes soliton solution $A$ of the effective NLS equation. We choose the focusing case ($\sigma <0$) and selected the carrier frequency of the wavepacket in the interior
of one of the higher energy bands far away from band gaps.

We select the peirodic potentials
\beq\label{E:Vsig_ex}
V(x)=\cos(x_1)\cos(x_2), \ \sigma(x) = \cos(x_1)\cos(x_2)-2
\eeq
and define the carrier Bloch wave by choosing
\beq\label{E:kn_ex}
k_0=(0.4,0)^T, \ n_0=4,
\eeq
which results in $\omega_0 \approx 2.075, v_g \approx (2.5083,0)^T$ and in the effective NLS coefficients $D^2\omega_{n_0}(k_0)\approx 1.5854*I_{2\times 2}$ and $\nu \approx 0.04905.$

\begin{figure}[ht!]
\includegraphics[scale=0.6]{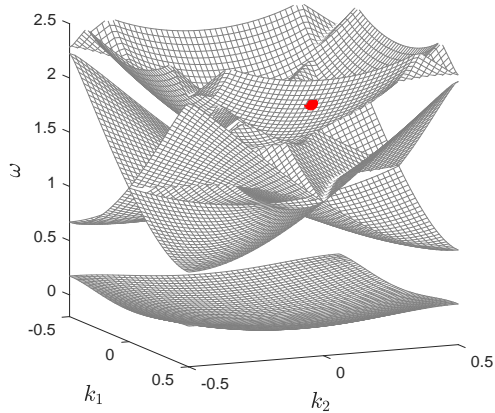}
\caption{The band structure for $V(x)=\cos(x_1)\cos(x_2)$ and the marked point $(k_0,\omega_{n_0}(k_0))$ with $k_0=(0.4,0)^T$ and $n_0=4$.}
\label{Fig:band_str}
\end{figure}

The band structure corresponding to the potential $V$ in \eqref{E:Vsig_ex} is plotted in Fig. \ref{Fig:band_str}. It was computed using a fourth order centered finite difference discretization of the eigenvalue problem \eqref{E:Bloch-gen}.

\begin{figure}[ht!]
\includegraphics[scale=0.8]{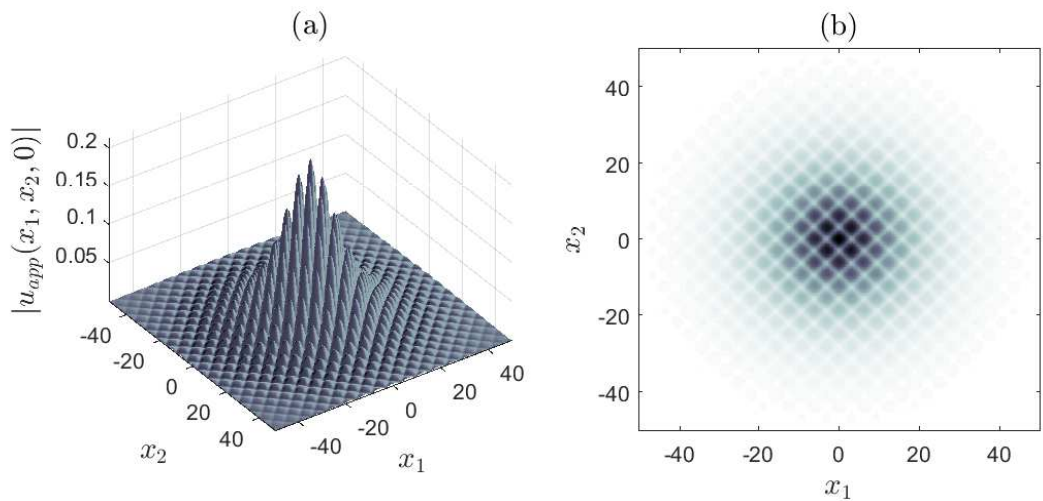}
\caption{The modulus of the asymptotic approximation $u_\text{app}(x,t)$ at $t=0$ with $\eps=0.1$. Side view in (a) and top view in (b).}
\label{Fig:as_prof_eps_p1}
\end{figure}

In Fig. \ref{Fig:as_prof_eps_p1} we plot the asymptotic approximation $u_\text{app}(x,0)$ from \eqref{E:uapp} with $\eps=0.1$, where $A(X,T)=e^{\ri T}R(|X|)$ is the Townes soliton of the effective NLS \eqref{E:NLS}. The radial profile $R$ is a solution of a ordinary differential equation in the radius $r:=|x|$. It is computed using a shooting method and the 4-5th oder Runge-Kutta Matlab scheme ode45.

Next, we solve the GP \eqref{E:GP} with the initial data $u(x,0)=u_\text{app}(x,0)$. We employ here the second order split step method, so called Strang splitting \cite{WH86}. In this method the problem is split into the linear constant coefficient part $\ri \pa_t u + \Delta u=0$ and the rest $\ri \pa_t u -V(x)u-\sigma(x)|u|^2u=0$. The former is solved in Fourier variables $\widehat{u}(k,t):=\frac{1}{(2\pi)^2}\int_{\R^2}u(x,t)e^{-\ri k\cdot x}\dd x$ exactly via $\uhat(k,t) = e^{-\ri |k|^2t}\uhat(k,0)$. In the implementation the Fourier transform is, of course, replaced by the fast Fourier discrete transform applied to the values of $u(\cdot,t)$ on a discrete grid. The latter part can be solved exactly in physical space via $u(x,t)=e^{-\ri (V(x)+\sigma(x)|u(x,0)|^2)t}u(x,0)$. These two problems are then suitable combined for each temporal discretization interval $t\in [n dt, (n+1) dt]$.

For the temporal discretization we use $dt =0.02$ and in space we discretize with $dx_1=dx_2 \approx 0.203$. The computational box is selected relatively large, namely $[-20\pi-\tfrac{5}{4\eps^2},20\pi+\tfrac{5}{4\eps^2}]\times [-40\pi,40\pi]$. This guarantees that the pulse remains well inside the box within the computaional time interval $[0,\tfrac{1}{\eps^2}]$.

\begin{figure}[ht!]
\includegraphics[scale=0.8]{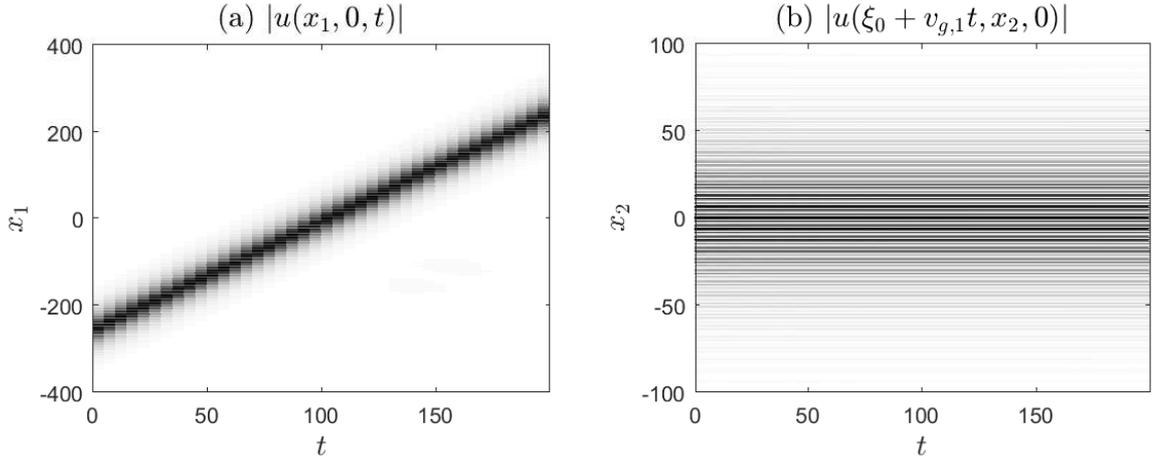}
\caption{Time evolution of the solution modulus $|u|$ for the example in \eqref{E:Vsig_ex} and \eqref{E:kn_ex} with $\eps =0.05$. The initial data are centered at $(x_1,x_2)=(\xi_0,0)$ with $v_g \approx (2.5083,0)$. (a) $|u(x_1,0,t)|$; (b) $|u(\xi_0+v_{g,1}t,x_2,t)|$. The apparent discontinuities in (a) are only due to the sampling of the numerical data for plotting.}
\label{Fig:waterfall_eps_p05}
\end{figure}

For $\eps=0.05$ the approximation via $\uapp(x,t)$ is very good and the resulting solution is close to a solitary wave. The temporal evolution is plotted in Fig. \ref{Fig:waterfall_eps_p05}. In (a) we plot the solution at $x_2=0$ and in (b) at $x_1=\xi_0+v_{g,1}t$, where $\xi_0$ is the $x_1$-position of the maximum of the initial pulse $|u(x,0)|$, i.e. $\xi_0+v_{g,1}t$ is the $x_1$-position of the maximum of $|u_\text{app}(\cdot,t)|$.

We also study the $\eps$-convergence of the error in the supremum norm $\|u(\cdot,\eps^2t)-u_\text{app}(\cdot,\eps^2t)\|_{C^0_b}$. The error is computed for the five values $\eps=0.3, 0.2, 0.1, 0.05, 0.03$. The convergence is $c\eps^{2.25}$, i.e. slightly better than predicted in Theorem \ref{T:main}.
\begin{figure}[ht!]
\includegraphics[scale=0.6]{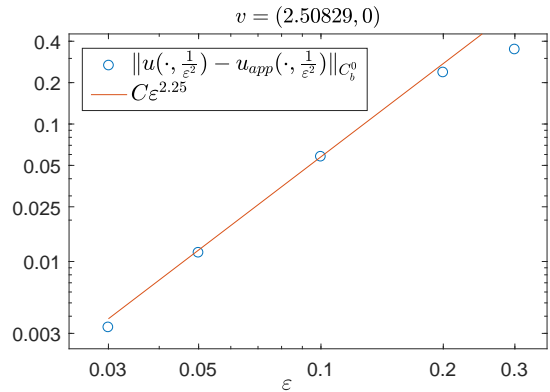}
\caption{$\eps$-convergence of the error in the supremum norm.}
\label{Fig:eps_conv_vg_2p5}
\end{figure}

\begin{sloppypar}
Movies of the time evolution for the $\eps$-values $\eps=0.1$ and $\eps =0.05$ are available at \url{http://www.mathematik.uni-dortmund.de/~tdohnal/2DPNLS-pulse-movies.html}. At $\eps=0.1$ some initial shedding of radiation can be observed, after which the pulse seems to become close to a solitary wave. At $\eps=0.05$ almost no radiation is visible ``with the naked eye''.
\end{sloppypar}

\section*{Acknowledgments}
This research is supported by the \emph{German Research Foundation}, DFG grant No. DO1467/3-1.
D.R. is supported by the SFB/TRR 191 `Symplectic Structures in Geometry,
Algebra and Dynamics', funded by the DFG.

\bibliographystyle{plain}
\bibliography{biblio-PNLS-dD}

\end{document}